\documentclass[1pt]{article}
\usepackage{graphicx}%
\usepackage{multirow}%
\usepackage{amsmath,amssymb,amsfonts}%
\usepackage{amsthm}%
\usepackage{mathrsfs}%
\usepackage[title]{appendix}%
\usepackage{xcolor}%
\usepackage[export]{adjustbox}
\usepackage{bm}
\usepackage{footnotehyper}
\usepackage{amssymb,amsmath,amsthm}
\usepackage{amscd}
\usepackage{hyperref}
\usepackage{lipsum}
\usepackage{epsfig}
\usepackage{amsfonts}
\usepackage{amssymb}
\usepackage{amsmath,enumerate}
\usepackage{commath}
\usepackage{euscript}
\usepackage{enumitem}
\usepackage{amscd}
\usepackage{xcolor}
\usepackage{makecell}
\usepackage{hhline}
\usepackage[numbers,sort&compress]{natbib}
\usepackage{lscape}
\usepackage{textcomp}%
\usepackage{manyfoot}%
\usepackage{booktabs}%
\usepackage{algorithm}%
\usepackage{algorithmicx}%
\usepackage{algpseudocode}%
\usepackage{listings}%
\usepackage{multirow}
\usepackage{amsmath, amsthm, amssymb}
\usepackage{graphicx}
\usepackage{algorithm,setspace}
\usepackage{rotating}
\usepackage{enumitem}
\usepackage{algpseudocode}
\newtheorem{corollary}{Corollary}[section]  
\usepackage{commath}
\usepackage{booktabs} 

\newtheorem{lemma}{Lemma}  

\usepackage{hyperref}
\hypersetup{
	colorlinks=true,
	linkcolor=blue,
	filecolor=blue,
	urlcolor=blue,
}
\theoremstyle{thmstyleone}%
\newtheorem{theorem}{Theorem}
\theoremstyle{thmstyletwo}%
\newtheorem{remark}{Remark}%

\theoremstyle{thmstylethree}%

\begin{document}
	\title{Estimation of scale parameters in scale mixtures of exponential distributions}
	\author{Somnath Mondal\footnote
		{\baselineskip=10pt
			~somnathmondal@iitbhilai.ac.in; somnathmondal1498@gmail.com}~ 
		and Lakshmi Kanta Patra\\
		Department of Mathematics\\
		Indian Institute of Technology Bhilai, Bhilai, India-491002}
	
	\date{}
	\maketitle
	\begin{abstract}
		The scale mixture of the exponential distribution provides a flexible framework for modelling lifetime and reliability data. This model is widely used in survival analysis, biomedical studies, statistical finance, and other related disciplines. In this work, we investigate the estimation of the ordered scale parameter of two scale mixture of exponential distributions under Stein loss and symmetric loss functions. Under certain conditions, we prove the inadmissibility of the affine equivariant estimator and exhibit several improved estimators. Consequently, we propose a class of estimators that uniformly dominate the best affine equivariant estimators (BAEE). Furthermore, we have proved that the boundary estimator of this class is a generalized Bayes estimator. As an application, we have proposed improved estimators for the ordered scale parameters of the multivariate Lomax and exponential inverse Gaussian distributions. For each case, we have conducted a simulation study to compare the risk performance of the improved estimators. Finally, we have given two real-life data analysis for implementation purposes.\\
		
	\end{abstract}
	
	
	\section{Introduction} 
	The scale mixture of exponential distributions has attracted considerable attention because of its wide range of applications in reliability theory, survival analysis, actuarial science, and lifetime data analysis. This distribution also includes the multivariate Lomax distribution and the exponential-inverse Gaussian (E–IG) distribution; see \cite{lindley1986multivariate}, \cite{nayak1987multivariate}, and \cite{petropoulos2006estimation}.  Let $X_1,X_2,\dots,X_{p}$ have the joint density 
	\begin{equation}\label{mix1}
		f(x_1,x_2,\dots,x_{p};\mu,\sigma)=\int_{0}^{\infty}\frac{\tau^{p}}{\sigma^{p}}e^{-\frac{\tau}{\sigma}\sum \limits_{i=1}^{p}(x_i-\mu)}I_{(\mu,\infty)}(x_{(1)})dH(\tau)
	\end{equation}
	where $x_{(1)}=\min\left\{x_1,x_2,\dots,x_p\right\}$, $I_{(a,b)}(.)$ is the indicator function and $\mu \in \mathbb{R}$, $\sigma>0$ are unknown. Here $\tau$ is a mixing random variable  whose  distribution function is $H(\cdot)$. The interpretation of
	the equation (\ref{mix1}) is, given $\tau>0$,  $X_1,X_2,\dots,X_{p}$ are i.i.d. exponential random variables with location parameter $\mu$ and scale parameter $\sigma/\tau$ that is $Exp(\mu,\sigma/\tau)$. This model is named a scale mixture of exponential distributions. It was originally introduced by \cite{lindley1986multivariate} to assess the reliability of parallel and series systems in terms of their component reliabilities. The authors have argued that the lifetimes of $p$ independent components will not be independent when they operate in the same laboratory environment. The common laboratory environment may influence the lifetime of the components. To address this issue, \cite{lindley1986multivariate} proposes a change of scale parameter $\sigma$ to $\sigma/\tau$. If  $\tau$ is degenerate at $1$ then we get $X_1,\dots,X_p$ are i.i.d. exponential distribution $Exp(\mu,\sigma)$. Recently, some work has been done on improved estimation of the scale parameter and the function scale parameter of a scale mixture of exponential distributions; see \cite{petropoulos2005estimation}, \cite{petropoulos2010class}, and \cite{petropoulos2020improved}. In this work, we study component-wise estimation of the ordered scale parameter of the two scale mixture exponential distribution.  
	
	We will now briefly discuss the literature on estimating ordered location and scale parameters of various probability distributions..
	For a review of the estimation of restricted parameters and its application, we refer to \cite{silvapulle2005constrained} and \cite{van2006restricted}. 
	Improved estimation of the ordered location parameter of a two-exponential distribution under a squared error loss function has been investigated by \cite{pal1992order}. They have proved the inadmissibility of standard estimators. \cite{misra1994estimation} investigated component-wise estimation of the ordered location parameter of two exponential distributions with known scale parameters. They have shown that unrestricted minimum risk estimators are inadmissible. \cite{vijayasree1995componentwise} studied the component-wise estimation of ordered scale and location parameters of several exponential distributions under a quadratic loss function. By proposing improved estimators, they showed that the unrestricted minimum risk equivariant estimator is inadmissible. They also proved that restricted MLE is inadmissible. Smooth improved estimators of the ordered scale parameter of the two gamma distributions are investigated by \cite{misra2002smooth} under squared error loss function. They proved that the proposed estimators are generalised Bayes. \cite{chang2008estimation} discussed the estimation of the linear function of ordered scale parameters of two gamma distributions under the entropy loss function. They have given a necessary and sufficient condition under which MLE dominates the crude unbiased estimator. \cite{jana2015estimation} dealt with estimating the ordered scale parameters of two exponential distributions with a common location parameter $\mu$. The authors obtained the UMVUE and showed that the restricted MLE dominates the usual MLE. Additionally, the authors proved an inadmissibility result as an application of the inadmissibility result they obtained for classes of scale and affine equivariant estimators. \cite{petropoulos2017estimation} investigated the estimation of ordered scale parameters of two multivariate Lomax distributions with unknown locations under a quadratic loss function. By proposing various improved estimators, the author showed that the best equivariant estimators are not admissible. \cite{patra2021componentwise} studied the component-wise estimation of the ordered scale parameter of two exponential distributions under a general scale-invariant loss function. They have proved the inadmissibility of the best affine equivariant estimators and restricted MLE by proposing several dominating estimators. For some more recent works on estimation of ordered parameters of exponential distributions, we refer to \cite{bobotas2019estimation}, \cite{patra2021componentwise}, \cite{kayal2024estimating}, \cite{jena2024estimating}, and \cite{bajpai2025improved}.     
	
	In the literature, much attention has not been given to the estimation of the ordered scale parameter in the scale-mixture of exponential distributions. This manuscript will investigate component-wise estimation of the ordered scale parameter of a two-scale mixture of exponential distributions. We consider the model: for a given $\tau>0$, let $X_1,X_2,\dots,X_{p_1}$ and $Y_1,Y_2,\dots,Y_{p_2}$ be a random sample taken from the exponential distribution $Exp\left(\mu_1,\frac{\sigma_1}{\tau}\right)$ and $Exp\left(\mu_2,\frac{\sigma_2}{\tau}\right)$ respectively where $-\infty<\mu_1,~\mu_2<\infty$ and $0<\sigma_1\leq\sigma_2$ are unknown parameters. 
	The complete and sufficient statistic is $(S_1,X,S_2,Y)$ (see \cite{lehmann2006theory}), where $S_1=\sum_{i=1}^{p_1}(X_i-X_{(1)})$, $X\equiv X_{(1)}=\min\left(X_1,X_2,\dots,X_{p_1}\right)$, $S_2=\sum_{i=1}^{p_2}(Y_i-Y_{(1)})$ and $Y\equiv Y_{(1)}=\min\left(Y_1,Y_2,\dots,Y_{p_2}\right)$. Given $\tau >0$, the statistics $S_1$, $S_2$, $X$ and $Y$ are independent random variables distributed as,
	\begin{align}
		& S_1|\tau \sim Gamma(p_1-1,\sigma_1/\tau), \quad X|\tau \sim \mathcal{E}\left(\mu_1,\sigma_1/(p_1\tau)\right) \notag \\
		& S_2|\tau \sim Gamma(p_2-1,\sigma_2/\tau), \quad Y|\tau \sim \mathcal{E}\left(\mu_2,\sigma_2/(p_2\tau)\right)
	\end{align}
	We consider the estimation of $\sigma_1$ and $\sigma_2$ with respect to the loss function 
	\begin{align*}
		\mbox{Symmetric loss:~}& L_1(\delta_i,\sigma_i)= \frac{\delta_i}{\sigma_i} +\frac{\sigma_i}{\delta_i} -2\\
		\mbox{Stein’s loss:~}& L_2(\delta_i,\sigma_i)= \frac{\delta_i}{\sigma_i} - \ln \left(\frac{\delta_i}{\sigma_i}\right)-1
	\end{align*}
	where $\delta_i$ is an estimator of $\sigma_i, i=1,2$. We invoke the principle of invariance and consider the affine group of transformations
	$G_{a_i,b_i}=\{g_{a_i,b_i}(x)=a_ix+b_{i}, j=1,\dots,p_i\}, ~i=1,2$. After some simplification, the form of the affine equivariant estimator is obtained as $cS_i$, where $c$ is a positive constant. 
	\begin{lemma}
		\begin{enumerate}
			\item [(i)] Under the $L_1(\cdot)$ loss function the BAEE of $\sigma_i$ is $\delta_{{1i}}=c_{i}S_i$, where $c_{i}=\sqrt{\frac{E(\tau)}{(p_i-1)(p_i-2)E(1/\tau)}}$ with $p_i\geq 3$ for $i=1, 2$.
			
			\item [(ii)] Under the $L_2(\cdot)$ loss function and for $i=1, 2$, the BAEE of $\sigma_i$ is $\delta_{{2i}}=d_{i}S_i$ with $d_{i}=\frac{1}{(p_i-1)E(1/\tau)}$ provided $p_i\geq2$.
		\end{enumerate}
	\end{lemma}
	
	The goal of this study to derive estimators of $\sigma_i$ for $i=1,2$ that dominate the BAEE under the order restriction $\sigma_1 \le \sigma_2$. The main contribution of this article is as follows.
	\begin{itemize}
		\item[(i)] We obtain BAEE of $\sigma_i,~i=1,2$ with respect to the loss function $L_1(\cdot)$  and $L_2(\cdot)$. Our goal is to propose estimators that dominate BAEE when $\sigma_1 \le \sigma_2$. We have derived several estimators whose risk is uniformly smaller than BAEE. Further, we have obtained a class of improved estimators and shown that the boundary estimator of this class is generalized Bayes. 
		
		\item[(ii)] As an application, we study the estimation of ordered scale parameters of the multivariate Lomax distribution and the exponential inverse Gaussian distribution.
		\item[(iii)] To compare the risk performance of the proposed estimators, we have conducted a detailed simulation study. The relative risk improvements of the proposed estimators over BAEE have been tabulated. 
	\end{itemize}

	The rest of the paper is organized as follows. In Section \ref{se2} of this paper, we have proposed improved estimators of the scale parameter $\sigma_1$ under two scale invariant loss functions $L_1(\cdot)$ and $L_2(\cdot)$ which dominate the BAEE. Also, we have obtained a class of estimators that improve upon the BAEE. In Section \ref{se3}, the improved estimation of $\sigma_2$ is considered. As an application, we have derived the improved estimators for scale parameters of the multivariate Lomax distribution (see in Subsection \ref{se4.1}) and exponential inverse Gaussian (E-IG) distribution (see Subsection \ref{se4.2}). Also, a simulation study is conducted to compare the performance of the proposed estimators. Finally, we have provided two real life data analysis in Section \ref{sec5}. 
	
	\section{Improved estimation for $\sigma_1$}\label{se2}
	Here we will consider the improved estimation of $\sigma_1$ when $\sigma_1\le \sigma_2$. We consider a class of estimators of the form 
	\begin{equation}\label{W}
		\mathcal{D}_1=\left\{\delta_{\varphi_1} = \varphi_1\left(W\right)S_1; \\ W=\frac{S_2}{S_1},\  \mbox{$\varphi_1(\cdot)$ is a positive function} \right\}.
	\end{equation} 
	Observe that the BAEE lies in this class. In the next theorem, we find an improved estimator that dominates the $\delta_{\varphi_1}$. As an application of this theorem, we derive an estimator that dominates the BAEE $\delta_{j1}$ for $j=1,2.$ 
	\begin{theorem}\label{Th1}
		\begin{enumerate}
			\item [(i)] Risk of the estimator $\delta^1_{1S1}(X,S)=\min\left\{\varphi_1(W),\varphi_{11}(W)\right\}S_1$
			is nowhere larger than that of $\delta_{\varphi_1}$ under symmetric loss function $L_1(\cdot)$ provided $P\left(\varphi_{11}(W)<\varphi_1(W)\right)>0$, where $\varphi_{11}(W)=\frac{(1+W)\sqrt{E(\tau)}}{\sqrt{(p_1+p_2-2)(p_1+p_2-3)E(1/\tau)}}$ with $p_1+p_2>3$.
			\item [(ii)] Risk of the estimator $\delta^2_{1S1}(X,S) = \min\left\{\varphi_1(W),\varphi_{12}(W)\right\}S_1$
			is nowhere larger than that of $\delta_{\varphi_1}$ under the Stein loss function $L_2(\cdot)$ provided $P\left(\varphi_{12}(W)<\varphi_1(W)\right)>0$, where $\varphi_{12}(W)=\frac{(1+W)}{(p_1+p_2-2)E(\frac{1}{\tau})}$ with $p_1+p_2>2$. 
		\end{enumerate}
	\end{theorem}
	\noindent\textbf{Proof.}~
	\textbf{(i)} We can easily seen that the risk function of the estimator is of the form (\ref{W}) depend on the unknown parameter through $\eta=\frac{\sigma_1}{\sigma_2} \le 1$ and it can be written as 
	\begin{align*}
		R(\delta_{\varphi_1};\eta)
		&=E\left\{E\left[(\varphi_1(W)S_1/\sigma_1) +\frac{1}{(\varphi_1(W)S_1/\sigma_1)} -2 \bigg| W\right]\right\}
	\end{align*}
	The inner conditional risk is minimized at 
	\begin{equation}\label{ph1}
		\varphi_{1}(w;\eta)=\left(\frac{E\left(1/V_1|W=w\right)}{E\left(V_1|W=w \right)}\right)^{\frac{1}{2}}.
	\end{equation}
	We can observe that the conditional distribution of $V_1=\frac{S_1}{\sigma_1}$, given $W=w>0$ and $T=\tau>0$ is proportional to 
	\begin{equation*}\label{den1}
		v_1^{p_1+p_2-3}e^{-\tau v_1(1+w\eta)}\tau^{p_1+p_2-2},\ w>0,\ v_1>0,
	\end{equation*}
	and we have 
	\begin{equation*}\label{phi2}
		\varphi_{1}(w;\eta)=\left(\frac{\int_{0}^{\infty}\int_{0}^{\infty}v_1^{p_1+p_2-4}e^{-\tau v_1(1+w\eta)}\tau^{p_1+p_2-2}dv_1dH(\tau)}{\int_{0}^{\infty}\int_{0}^{\infty}v_1^{p_1+p_2-2}e^{-\tau v_1(1+w\eta)}\tau^{p_1+p_2-2}dv_1dH(\tau)}\right)^{\frac{1}{2}}
	\end{equation*}
	Using the transformation $z_1=\tau v_1(1+w\eta)$, we obtain equation (\ref{phi2}) in the following form,
	\begin{align*}
		\varphi_{1}(w;\eta)= &  (1+w\eta)\left(\frac{\int_{0}^{\infty}\tau\int_{0}^{\infty}z_1^{p_1+p_2-4}e^{-z_1}dz_1dH(\tau)}{\int_{0}^{\infty}\frac{1}{\tau}\int_{0}^{\infty}z_1^{p_1+p_2-2}e^{-z_1}dz_1dH(\tau)}\right)^{\frac{1}{2}} =  (1+w\eta)\left(\frac{E(\tau)\left(E(1/\tau)\right)^{-1}}{(p_1+p_2-2)(p_1+p_2-3)}\right)^{\frac{1}{2}}\\
		\le &  (1+w)\left(\frac{E(\tau)\left(E(1/\tau)\right)^{-1}}{(p_1+p_2-2)(p_1+p_2-3)}\right)^{\frac{1}{2}}=\varphi_{11}(w).
	\end{align*}
	Now using the convexity of $R(\delta_{\varphi_1};\eta)$ we get the result. Proof of $(ii)$ similar to $(i)$. So we omit it. 
	\begin{corollary}
		Risk of the estimator $\delta_{11}^1=\min\left\{ c_1,\alpha_1 (1+W)\right\}S_1$ is nowhere larger than that of $\delta_{11}$ under the loss function $L_1(\cdot)$ provided $\alpha_1<c_1$, where $\alpha_1=\left(\frac{E(\tau)}{(p_1+p_2-2)(p_1+p_2-3)E(1/\tau)}\right)^{1/2}$, 	$p_1\geq 3$ and $p_2\geq 1$.
	\end{corollary}
	\begin{corollary}
		The estimator $\delta_{11}^2=\min\left\{ d_1,\beta_1 (1+W)\right\}S_1$ dominates the BAEE $\delta_{12}$ with respect to the loss function  $L_2(\cdot)$  provided $\beta_1<d_1$, where $\beta_1=\frac{1}{(p_1+p_2-2)E(1/\tau)}$, $p_1\geq 2$ and $p_2\geq 1$.
	\end{corollary}
	
	We next derive a class of improved estimators. For this purpose, we applied the integral expression of risk difference (IERD) approach of \cite{kubokawa1994double} and \cite{kubokawa1994unified}. In the following theorem, we give sufficient conditions under which we obtain a class of improved estimators. 
	\begin{theorem}\label{thIERD}
		\begin{enumerate}
			\item [(i)] Suppose $p_1\geq 3$ and $p_2\geq 1$, and the function $\varphi_1(w)$ satisfies the following conditions:
			\begin{enumerate}
				\item [(a)] $\varphi_1(w)$ is non-decreasing in $w$ and $\lim\limits_{w\rightarrow\infty}\varphi_1(w)=\left(\frac{E(\tau)}{(p_1-1)(p_1-2)E(1/\tau)}\right)^{1/2}$
				\item [(b)] $\varphi_1(w)\geq\varphi_{*}^1(w)$, where $\varphi_{*}^1(w)= \left(\frac{E(\tau) \int_{0}^{\frac{w}{1+w}} s^{p_2-2}(1-s)^{p_1-3} ds} {E(1/\tau)(p_1+p_2-2)(p_1+p_2-3) \int_{0}^{\frac{w}{1+w}} s^{p_2-2}(1-s)^{p_1-1} ds}\right)^{1/2}$.
			\end{enumerate}
			Then the estimator $\delta_{\varphi_1}$ (defined in (\ref{W})) dominates $\delta_{11}$ under the symmetric loss function $L_1(\cdot)$. 
			\item [(ii)] Let $p_1\geq 2$ and $p_2\geq 1$, and the function $\varphi_1(w)$ satisfies the following conditions:
			\begin{enumerate}
				\item [(a)] $\varphi_1(w)$ is non-decreasing in $w$ and $\lim\limits_{w\rightarrow\infty}\varphi_1(w)=\frac{1}{(p_1-1)E\left(1/\tau\right)}$
				\item [(b)] $\varphi_1(w)\geq\varphi_{*}^2(w) = \frac{ \int_{0}^{\frac{w}{1+w}} s^{p_2-2}(1-s)^{p_1-2}ds} {E(1/\tau)(p_1+p_2-2) \int_{0}^{\frac{w}{1+w}} s^{p_2-2}(1-s)^{p_1-1}ds}$.
			\end{enumerate}		
			Then the risk of $\delta_{\varphi_1}$ is nowhere larger than that of $\delta_{21}$ with respect to the Stein loss function $L_2(\cdot)$.
		\end{enumerate}
	\end{theorem}
	\noindent \textbf{Proof.} The proof is deferred to Appendix \ref{proof1}.

	\begin{remark}
		The boundary estimator $\delta_{\varphi_{*}^1}$ and  $\delta_{\varphi_{*}^2}$ are Brewster- Zidek-type (\cite{brewster1974improving}) estimator for $\sigma_1$ under the loss functions $L_1(\cdot)$ and $L_2(\cdot)$ respectively. 
		
	\end{remark}
	\begin{remark}
		Now we prove that $\delta_{\varphi_{*}^1}$ is a generalized Bayes estimator of $\sigma_1$ with respect to $L_1(\cdot)$. We consider the prior distribution as 
		\begin{align*}
			\pi(\sigma_1,\sigma_2,\mu_1,\mu_2)=\frac{1}{\sigma_1\sigma_2},\quad \sigma_1\leq\sigma_2,\ \mu_1,\mu_2 \in \mathbb{R}.
		\end{align*} 
		The corresponding posterior distribution, for given $\tau> 0$, is
		\begin{equation}\label{post}
			\pi_1(\sigma_1,\sigma_2,\mu_1,\mu_2\big| X,S_1,Y,S_2)\propto \frac{\tau^{p_1+p_2-2}}{\sigma_1^{p_1}\sigma_2^{p_2}}e^{-\frac{\tau S_1}{\sigma_1}-\frac{\tau S_2}{\sigma_2}} \frac{p_1\tau}{\sigma_1} e^{-\frac{p_1\tau}{\sigma_1}(X-\mu_1)} \frac{p_2\tau}{\sigma_2}e^{-\frac{p_2\tau}{\sigma_2}(Y-\mu_2)},
		\end{equation}
		where $ \mu_1\leq x,\ \mu_2\leq y,\ 0< \sigma_1 \leq \sigma_2$. For the symmetric loss function $L_1(\cdot)$, the generalized Bayes estimator of $\sigma_1$ is obtained as follows
		$$\delta_{1B}^1= \left(\frac{E\left(\sigma_1 \big| X,S_1,Y,S_2\right)}{E\left(\frac{1}{\sigma_1} \big|  X,S_1,Y,S_2\right)}\right)^{1/2}.$$
		After some calculation, it is found that the generalized Bayes estimator coincides with the estimator $\delta_{\varphi_{*}^1}$. 		
	\end{remark}
	\begin{remark} By using a similar argument, we can prove that $\delta_{\varphi_{*}^2}$ is generalized Bayes estimator of $\sigma_1$ under the Stein loss function $L_2(\cdot)$, with respect to the same prior distribution $\pi(\sigma_1,\sigma_2,\mu_1,\mu_2)$.
	\end{remark}
	\noindent Now we  define a bigger class of estimators based on the statistics $(S_1,S_2, X)$  of the form 
	\begin{equation}\label{WW1}
		\mathcal{D}_2= \left\{ \delta_{\varphi_2}=\varphi_2\left(W,W_1\right)S_1;\ W=\frac{S_2}{S_1},\  W_1=\frac{X}{S_1} \right\}.
	\end{equation}
	In the next theorem we will propose an estimator that dominates the estimator $\delta_{\varphi_2}$. As a consequences of this theorem we will obtain an estimator that dominates BAEE of $\sigma_1$. 
	\begin{theorem}\label{ths21}
		\begin{enumerate}
			\item [(i)] Risk of the estimator
			\begin{equation}
				\delta_{1S2}^{1}
				= \left\{
				\begin{array}{ll}
					\displaystyle
					\min\left\{\varphi_2(W,W_1),\,\varphi_{21}(W,W_1)\right\}S_1,
					& \text{if } W_1>0, \\[2ex]
					\displaystyle
					\varphi_2(W,W_1)S_1,
					& \text{otherwise.}
				\end{array}
				\right.
			\end{equation}
			is nowhere larger than that of $\delta_{\varphi_2}$ under the loss function $L_1(\cdot)$ provide $P(\varphi_2(W,W_1)$ $>\varphi_{21}(W,W_1))>0$, where $\varphi_{21}(W,W_1)=\frac{(1 + W + p_1 W_1)}{\sqrt{(p_1 + p_2 - 1)(p_1 + p_2 - 2)}} \min\left\{ 
			\left( \frac{E(\tau^{p_1 + p_2 + 1})}
			{E(\tau^{p_1 + p_2 - 1})} \right)^{1/2}, \left(\frac{E(\tau)}	{E(\tau^{-1})} \right)^{1/2}
			\right\}$ and $p_1+p_2\geq 3$.
			
			\item [(ii)] Risk of the estimator
			\begin{equation}
				\delta_{1S2}^{2} = \begin{cases} 
					\min\left\{\varphi_2W,W_1),\varphi_{22}(W,W_1)\right\}S_1, & \ W_1>0\\
					\varphi_2(W,W_1)S_1, & \text{otherwise}
				\end{cases}
			\end{equation} is nowhere larger than that of $\delta_{\varphi_2}$ with respect to the loss function $L_2(\cdot)$ provide $P(\varphi_2(W,W_1)>	\varphi^2_2(W,W_1))>0$, where $	\varphi_{22}(W,W_1)=\frac{(1 + W + p_1 W_1)}{p_1 + p_2 - 1} \min\Bigg\{\frac{E(\tau^{p_1 + p_2})}{E(\tau^{p_1 + p_2 - 1})} , \frac{1}{E(\tau)} \Bigg\}$ and $p_1+p_2\geq 2$.
		\end{enumerate}
	\end{theorem}
	\noindent\textbf{Proof.} The proof is deferred to Appendix \ref{proof2}.
	\begin{corollary}
		The estimator $\delta_{12}^1=\min\left\{ c_1,\alpha_2 (1+W+p_1W_1)\right\}S_1$ dominates  $\delta_{11}$ with respect to the loss function $L_1(\cdot)$ provided $\alpha_2<c_1$, where $\alpha_2=\left(\frac{E(\tau^{p_1+p_2+1})}{(p_1+p_2-1)(p_1+p_2-2)E(1/\tau^{p_1+p_2-1})}\right)^{1/2}$,  $p_1\geq 3$ and $p_2\geq 1$.
	\end{corollary}
	\begin{corollary}
		Risk of the estimator $\delta_{12}^2=\min\left\{ d_1,\beta_2 (1+W+p_1W_1)\right\}S_1$ is nowhere larger than that of $\delta_{12}$ under the loss function $L_2(\cdot)$ provided $\beta_2<d_1$, where $\beta_2=\frac{E(\tau^{p_1+p_2})}{(p_1+p_2-1)E(1/\tau^{p_1+p_2-1})}$, $p_1\geq 2$ and $p_2\geq 1$.
	\end{corollary}
	Now the Theorem \ref{ths21} can be extended by using the information contained in the statistics $Y$, substituting $Y$ for $X$ within the class $\mathcal{D}_2$ in (\ref{WW1}) and define $W_2=\frac{Y}{S_1}$.
	\begin{theorem}\label{th2.4}
		\begin{enumerate}
			\item [(i)] Risk of the estimator
			\begin{equation}
				\delta_{1S3}^{1} = \begin{cases} 
					\min\left\{\varphi_3(W,W_2),\varphi_{31}(W,W_2)\right\}S_1, & \ W_2>0\\
					\varphi_3(W,W_2)S_1, & \text{otherwise}
				\end{cases}
			\end{equation}  
			is nowhere larger than that of $\delta_{\varphi_3}$ under the loss function $L_1(\cdot)$, provide $P(\varphi_3(W,W_2)> \varphi_{31}(W,W_2))>0$, where  $	\varphi_{31}(W,W_2)=\frac{(1 + W + p_2 W_2)}{\sqrt{(p_1 + p_2 - 1)(p_1 + p_2 - 2)}} \min\Bigg\{ \left( \frac{E(\tau^{p_1 + p_2 + 1})}
			{E(\tau^{p_1 + p_2 - 1})} \right)^{1/2}, \left(\frac{E(\tau)}	{E(\tau^{-1})} \right)^{1/2}
			\Bigg\}$ and $p_1+p_2> 2$.
			\item [(ii)]  Risk of the estimator
			\begin{equation}
				\delta_{1S3}^{2} = \begin{cases} 
					\min\left\{\varphi_3(W,W_2),\varphi_{32}(W,W_2)\right\}S_1, & \ W_2>0\\
					\varphi_3(W,W_2)S_1, & \text{otherwise}
				\end{cases}
			\end{equation} is nowhere larger than that of $\delta_{\varphi_3}$ under the loss function $L_2(\cdot)$ provide $P(\varphi_3(W,W_2)>\varphi_{32}(W,W_2))>0$, where $	\varphi_{32}(W,W_2)=\frac{(1 + W + p_2 W_2)}{p_1 + p_2 - 1} \min\Bigg\{ \frac{E(\tau^{p_1 + p_2})} {E(\tau^{p_1 + p_2 - 1})} , \frac{1}	{E(\tau^{-1})} 
			\Bigg\}$ and $p_1+p_2\geq 2$.
		\end{enumerate}
	\end{theorem}
	\begin{corollary}
		The estimator $\delta_{13}^1=\min\left\{ c_1,\alpha_3 (1+W+p_2W_2)\right\}S_1$ dominates $\delta_{11}$ with respect to the loss function $L_1(\cdot)$ provided $\alpha_3<c_1$, where $\alpha_3=\left(\frac{E(\tau^{p_1+p_2+1})}{(p_1+p_2-1)(p_1+p_2-2)E(1/\tau^{p_1+p_2-1})}\right)^{1/2}$, $p_1\geq 3$ and $p_2\geq 1$.
		
	\end{corollary}
	\begin{corollary}
		Risk of the estimator $\delta_{13}^2=\min\left\{ d_1,\beta_3 (1+W+p_2W_2)\right\}S_1$ is nowhere larger than that of $\delta_{12}$ under loss function $L_2(\cdot)$ provided $\beta_3<d_1$, where $\beta_3=\frac{E(\tau^{p_1+p_2})}{(p_1+p_2-1)E(1/\tau^{p_1+p_2-1})}$, $p_1\geq 2$ and $p_2\geq 1$. 
	\end{corollary}
	
	\section{Improved estimation for $\sigma_2$}\label{se3}
	In this section, we will derive estimators of the parameter $\sigma_2$ which will improve upon the BAEE fo $\sigma_2$ under the restriction $\sigma_1\leq\sigma_2$. We consider estimators of the form
	\begin{equation}\label{W*}
		\mathcal{C}_1= \left\{ \delta_{\psi_1} = \psi_1\left(U\right)S_2; \ U=\frac{S_1}{S_2},\  \mbox{$\psi_1(\cdot)$ is a positive function}\right\}
	\end{equation}
	and derive estimators for $\sigma_2$, which gives improvement over $\delta_{j2}$ for $j=1,2$.
	
	\begin{theorem}\label{tgh}
		\begin{enumerate}
			\item [(i)] Risk of the estimator
			$\delta_{2S1}^1(X,S) = 
			\max\left\{\psi_{1}(U),\psi_{11}(U)\right\}S_2,$
			is nowhere larger than that of $\delta_{\psi_1}$ under the loss function $L_1(\cdot)$ provided $P\left(\psi_{11}(U)>\psi_1(U)\right)>0$, where $\psi_{11}(U)=(1+U) \left(\frac{E(\tau)}{(p_1+p_2-2)(p_1+p_2-3)E(1/\tau)}\right)^{1/2}$ with $p_1+p_2>3$. 
			
			\item [(ii)] Risk of the estimator $\delta_{2S1}^2(X,S) = \max\left\{\psi_{1}(U),\psi_{12}(U)\right\}S_2,$
			is nowhere larger than that of $\delta_{\psi_1}$ under the loss function $L_2(\cdot)$ provided $P\left(\psi_{12}(U)>\psi_1(U)\right)>0$, where $\psi_{12}(U)=\frac{(1+U)}{(p_1+p_2-2)E(\frac{1}{\tau})}$ with $p_1+p_2\geq3$.
		\end{enumerate}	
	\end{theorem}
	\noindent\textbf{Proof.} Proof of this theorem is similar to the Theorem \ref{ths21}.
	\begin{corollary}
		The estimator $\delta_{21}^1=\max\left\{ c_2,\alpha_1^* (1+U)\right\}S_2$ dominates the BAEE $\delta_{21}$ under the loss function $L_1(\cdot)$ provided $\alpha_1^*<c_2$, where $\alpha_1^*=\left(\frac{E(\tau)}{(p_1+p_2-2)(p_1+p_2-3)E(1/\tau)}\right)^{1/2}$, $p_1\geq 1$ and $p_2\geq 3$.
	\end{corollary}
	\begin{corollary}
		Risk of the estimator $\delta_{21}^2=\max\left\{ d_2,\beta_1^* (1+U)\right\}S_2$ dominates the BAEE $\delta_{22}$ the loss function $L_2(\cdot)$ provided $\beta_1^*<d_2$, where $\beta_1^*=\frac{1}{(p_1+p_2-2)E(1/\tau)}$, $p_1\geq 1$ and $p_2\geq 2$.
	\end{corollary}
	
	As in the previous section, we derive an improved estimator for $\sigma_2$ within the class of estimators $\mathcal{C}_1$ using the IERD method of \cite{kubokawa1994double} and \cite{kubokawa1994unified}. We have the following theorem as follows
	\begin{theorem}\label{ThIERD}
		\begin{enumerate}
			\item [(i)] Let $p_1\geq 1$ and $p_2\geq 3$. Assume that the function $\psi_1(u)$ satisfies the following conditions:
			\begin{enumerate}
				\item [(a)] $\psi_1(u)$ is non-decreasing in $u$ and $\lim\limits_{u\rightarrow 0}\psi_1(u)=\left(\frac{E(\tau)}{(p_2-1)(p_2-2)E(1/\tau)}\right)^{1/2}$
				\item [(b)] $\psi_1(u)\leq\psi_{*}^1(u) = \left(\frac{E(\tau) \left[ \Gamma(p_2-2) - \Gamma(p_1+p_2-3) \int_{0}^{\frac{u}{1+u}} s^{p_1-2}(1-s)^{p_2-3} ds \right]} {E(1/\tau)  \left[\Gamma(p_2) - \Gamma(p_1+p_2-1) \int_{0}^{\frac{u}{1+u}} s^{p_1-2}(1-s)^{p_2-1} ds \right]}\right)^{1/2}$
				
			\end{enumerate}
			Then risk of the estimator $\delta_{\psi_1}$ is nowhere larger than that of $\delta_{12}$ under the loss function $L_1(\cdot)$.
			\item [(ii)] Suppose $p_1\geq 1$ and $p_2\geq 2$. Let $\psi_1(u)$ satisfies the following conditions:
			\begin{enumerate}
				\item [(a)] $\psi_1(u)$ is non-decreasing in $u$ and $\lim\limits_{u\rightarrow 0}\psi_1(u)=\frac{1}{(p_2-1)E\left(1/\tau\right)}$ 
				\item [(b)] $\psi_1(u)\leq\psi_{*}^2(u) = \frac{  \Gamma(p_2-1) - \Gamma(p_1+p_2-2) \int_{0}^{\frac{u}{1+u}} s^{p_1-2}(1-s)^{p_2-2} ds } {E(1/\tau) \left[\Gamma(p_2) - \Gamma(p_1+p_2-1) \int_{0}^{\frac{u}{1+u}} s^{p_1-2}(1-s)^{p_2-1} ds \right]}$.
				
			\end{enumerate}		
			Then risk of the estimator $\delta_{\psi_1}$ is nowhere larger than that of $\delta_{22}$ under the loss function $L_2(\cdot)$.
		\end{enumerate}
	\end{theorem}
	\noindent\textbf{Proof.} Proof of this theorem is similar to the Theorem \ref{thIERD}.

	\begin{remark}
		The boundary estimator $\delta_{\psi_{*}^1}$ and  $\delta_{\psi_{*}^2}$ are Brewster- Zidek-type (\cite{brewster1974improving}) estimator for $\sigma_2$ under the loss functions $L_1(\cdot)$ and $L_2(\cdot)$ respectively.
	\end{remark}
	\begin{remark}
		Now we prove that $\delta_{\psi_*^1}$ is a generalized Bayes estimator of $\sigma_2$ under the loss function $L_1(\cdot)$. We consider the prior distribution as
		\begin{align*}
			\pi(\sigma_1,\sigma_2,\mu_1,\mu_2)=\frac{1}{\sigma_1\sigma_2}, \quad 0<\sigma_1\leq\sigma_2, \ \mu_1,\mu_2 \in \mathbb{R}.
		\end{align*} 
		The corresponding posterior distribution, for given $\tau> 0$, is
		\begin{equation}\label{post1}
			\pi_2(\sigma_1,\sigma_2,\mu_1,\mu_2\big| X,S_1,Y,S_2)\propto \frac{\tau^{p_1+p_2-2}}{\sigma_1^{p_1}\sigma_2^{p_2}}e^{-\frac{\tau S_1}{\sigma_1}-\frac{\tau S_2}{\sigma_2}} \frac{p_1\tau}{\sigma_1} e^{-\frac{p_1\tau}{\sigma_1}(X-\mu_1)} \frac{p_2\tau}{\sigma_2}e^{-\frac{p_2\tau}{\sigma_2}(Y-\mu_2)},
		\end{equation}
		where $ \mu_1\leq x,\ \mu_2\leq y,\ 0< \sigma_1 \leq \sigma_2$. For the symmetric loss function $L_1(\cdot)$, the generalized Bayes estimator of $\sigma_2$ is obtained as follows
		$$\delta_{2B}^1= \left(\frac{E\left(\sigma_1 \big| X,S_1,Y,S_2\right)}{E\left(\frac{1}{\sigma_1} \big|  X,S_1,Y,S_2\right)}\right)^{1/2},$$
		After some calculation, we obtained that the generalized Bayes estimator coincides with the  estimator $\delta_{\psi_{*}^1}$. 
	\end{remark}
	\begin{remark} By using a similar argument, we can prove that $\delta_{\psi_{*}^2}$ is generalized Bayes estimator of $\sigma_2$ under the loss function $L_2(\cdot)$, with respect to the same prior distribution $\pi(\sigma_1,\sigma_2,\mu_1,\mu_2)$.
	\end{remark}
	
	We consider a class of estimator using the statistics $S_2$ and $Y$ as follows
	\begin{equation}
		\mathcal{C}_2= \left\{ \delta_{\psi_2} = \psi_2\left(U_1\right)S_2; \ U_1=\frac{Y}{S_2},\  \mbox{$\psi_2(\cdot)$ is a positive function} \right\}.
	\end{equation}
	\begin{theorem}\label{W2*}
		\begin{enumerate}
			\item [(i)] 
			Risk of the estimator
			\begin{equation*}
				\delta_{2S2}^{1} = \begin{cases} 
					\min\left\{\psi_2(U_1),\psi_{21}(U_1)\right\}S_2, & \ U_1>0\\
					\psi_2(U_1)S_2, & \text{otherwise}
				\end{cases}
			\end{equation*} is nowhere larger than that of $\delta_{\psi_2}$ under the loss function $L_1(\cdot)$ provide $P(\psi_2(U_1)$ $>\psi_{21}(U_1))>0$ where $\psi_{21}(U_1)=\frac{(1 + p_2 U_1)}{\sqrt{p_2(p_2 - 1)}} \min\Bigg\{ 
			\left( \frac{E(\tau^{ p_2 + 2})}
			{E(\tau^{p_2})} \right)^{1/2}, \left(\frac{E(\tau)}	{E(\tau^{-1})} \right)^{1/2}
			\Bigg\}$ and $p_2\geq 2$.
			\item [(ii)] Risk of the estimator 
			\begin{equation*}
				\delta_{2S2}^2 = \begin{cases} 
					\min\left\{\psi_2(U_1),\psi_{22}(U_1)\right\}S_2, & U_1>0,\\
					\psi_2(U_1)S_2, & \text{otherwise}
				\end{cases}
			\end{equation*}
			is nowhere larger than that of $\delta_{\psi_2}$ under the loss function $L_2(\cdot)$ provided $P(\psi_2(U_1)>\psi_{22}(U_1))>0$ where $\psi_{22}(U_1)=\frac{(1+p_2U_1)}{p_2E(\frac{1}{\tau})}$ and $p_2\geq 1$.
		\end{enumerate}
	\end{theorem}
	\noindent\textbf{Proof.} The theorem can be proved by using arguments similar to those in the proof of Theorem~\ref{ths21}.

	\begin{corollary}
		The estimator 
		\begin{align*}
			\delta_{22}^{1} = \begin{cases} 
				\min\left\{ c_2,\ \alpha_2^*(1+p_2U_1)\right\}S_2, & \ U_1>0\\
				c_2S_2, & \text{otherwise}
			\end{cases}
		\end{align*} dominates the BAEE $\delta_{21}$ under the loss function $L_1(\cdot)$ provided $\alpha_2^*<c_2$, where $\alpha_2^*=\left(\frac{E(\tau^{p_2+2})}{p_2(p_2-1)E(\tau^{p_2})}\right)^{1/2}$ and $p_2\geq 2$.
	\end{corollary}
	\begin{corollary}
		The estimator
		\begin{align*}
			\delta_{22}^{2} = \begin{cases} 
				\min\left\{ d_2,\ \beta_2^*(1+p_2U_1) \right\}S_2, & \ U_1>0\\
				d_2S_2, & \text{otherwise}
			\end{cases}
		\end{align*} dominates the BAEE $\delta_{22}$ under the loss function $L_2(\cdot)$ provided $\beta_2^*<d_2$, where $\beta_2^*=\frac{1}{p_2E(1/\tau)}$ and $p_2\geq 2$.
	\end{corollary}

	\section{Application}\label{se5}
	In this section, we will apply the results obtained in Sections \ref{se2} and \ref{se3} to two special lifetime distributions. 
	These distributions will be obtained by choosing a particular distribution of $\tau$. First, we assume that $\tau$ follows a Gamma distribution $\Gamma(b,1)$. Then the joint density of $X_1,X_2,\dots,X_{p}$ in (\ref{mix1}) is
	\begin{equation}\label{lom1}
		f_1(x_1,x_2,\dots,x_{p};\mu_1,\sigma_1)=\frac{\Gamma(p+b)}{\Gamma(b)\sigma_1^{p}}
		\frac{1}{\left(1+\frac{1}{\sigma_1}\sum_{i=1}^{p}(x_i-\mu_1)\right)^{p+b}}I_{(\mu_1,\infty)}(x_{(1)})
	\end{equation}
	which is the density function of  multivariate Lomax distribution. The Lomax distribution, also known as the Pareto Type II distribution. This distribution has been extensively used across various domains, including socioeconomics, reliability analysis, and life testing, biological and medical sciences, modeling business failure data etc. For details about the multivariate Lomax distribution we refer to  
	\cite{lindley1986multivariate}, \cite{holland2006traffic}, \cite{hassan2009optimum}. Secondly, when distribution of $\tau$ is an inverse Gaussian distribution $IG(m,n)$, then (\ref{mix1}) reduces to an exponential inverse Gaussian (E–IG) distribution having density function 
	\begin{equation}\label{ig1}
		\begin{split}
			f_2(x_1,x_2,\dots,x_{p};\mu_1,\sigma_1) =&\frac{\exp\left\{ \frac{n}{m}-n\sqrt{\frac{1}{m^2}+\frac{2}{\sigma_1n}\sum_{i=1}^{p}(x_i-\mu_1)}\right\}}{\sigma_1^{p}\left(\frac{1}{m^2}+\frac{2}{\sigma_1n}\sum_{i=1}^{p}(x_i-\mu_1)\right)^{p/2}}\times \\
			&\sum_{i=1}^{p-1}\frac{(p-1+i)!}{i!(p_1-1-i)!}\left[ 2n\sqrt{\left(\frac{1}{m^2}+\frac{2}{\sigma_1n}\sum_{i=1}^{p}(x_i-\mu_1)\right)}\right]^{-i}I_{(0,\infty)}(x_{(1)})
		\end{split}
	\end{equation}
	\cite{petropoulos2006estimation} studied quantile estimation for a scale mixture of exponential distributions with unknown location and scale parameters. As an application, the author obtained an improved estimator for the quantile of the exponential inverse Gaussian distribution. The E-IG distribution has been used in various domains, including reliability, actuarial science and survival analysis. Some work related to this distribution, we refer to \cite{bhattacharya1986ig}, \cite{hesselager1998exponential},  \cite{frangos2004modelling}.

	\subsection{Multivariate Lomax distribution}\label{se4.1}
	In this subsection, we will derive the improved estimator of the parameters $\sigma_1$ and $\sigma_2$ for multivariate Lomax distribution under the loss functions $L_1(\cdot)$ and $L_2(\cdot)$. In this case, we have $E(\tau)=b$ and $E(\tau^{-1})=\frac{1}{b-1}$ with $b\geq 2$. So that the BAEE of $\sigma_i$ under the loss function $L_1(\cdot)$ is $\delta_{{1i}}=c_iS_i$, where $c_i=\left(\frac{b(b-1)}{(p_i-1)(p_i-2)}\right)^{1/2}$ with $p_i\geq 3$ for $i=1,2$ and for the loss function $L_2(\cdot)$ we have the BAEE of $\sigma_i$ is $\delta_{{2i}}=d_iS_i$ with $d_i=\frac{b-1}{p_i-1}$ and $p_i\geq 2$.
	
	
	Now as an application of Theorem \ref{Th1}, \ref{ths21} and \ref{th2.4} the improved estimators of $\sigma_1$ are obtained for multivariate Lomax distribution in the following theorem.
	
	\begin{theorem}\label{4.1}
		\begin{enumerate}
			\item[(i)] Under $L_1(\cdot)$ loss function, we have $\varphi_{11}(W)=(1+W)\left(\frac{b(b-1)}{(p_1+p_2-2)(p_1+p_2-3)}\right)^{\frac{1}{2}}$,  $\varphi_{21}(W,W_1)=(1+W+p_1W_1)\left(\frac{b(b-1)}{(p_1+p_2-1)(p_1+p_2-2)}\right)^{1/2}$, $\varphi_{31}(W,W_2)=  \frac{(1+W+p_2W_2) \left(b(b-1)\right)^{1/2} }{\left((p_1+p_2-1)(p_1+p_2-2)\right)^{1/2}}$. The improved estimators of $\sigma_1$ are obtained as follows
			\begin{align*}
				& \delta^1_{11} = \min\left\{\varphi_{11}(W),\ c_{1}\right\}S_1,~~
				\delta_{12}^{1} = \begin{cases} 
					\min\left\{\varphi_{21}(W,W_1),c_{1}\right\}S_1, & \ W_1>0\\
					c_{1}S_1, & \text{otherwise}
				\end{cases}\\
				& 	\delta_{13}^{1} = \begin{cases} 
					\min\left\{\varphi_{31}(W,W_2),c_{1}\right\}S_1, & \ W_2>0\\
					c_{1}S_1, & \text{otherwise}
				\end{cases}
			\end{align*}
			
			\item [(ii)] Under the Stein loss function $L_2(\cdot)$ we have $\varphi_{12}(W)=(1+W)\frac{b-1}{(p_1+p_2-2)}$, $\varphi_{22}(W,W_1)=(1+W+p_1W_1)\frac{b-1}{p_1+p_2-1}$, $\varphi_{32}(W,W_1)=(1+W+p_2W_2)\frac{b-1}{p_1+p_2-1}$. The improved estimators of $\sigma_1$ are obtained as follows
			\begin{align*}
				& \delta^2_{11} =\min\left\{\varphi_{12}(W),d_{1}\right\}S_1,~~
				\delta_{12}^{2} = \begin{cases} 
					\min\left\{\varphi_{22}(W,W_1),d_{1}\right\}S_1, & \ W_1>0\\
					d_{1}S_1, & \text{otherwise}
				\end{cases}\\
				& 	\delta_{13}^{2} = \begin{cases} 
					\min\left\{\varphi_{32}(W,W_2),d_{1}\right\}S_1, & \ W_2>0\\
					d_{1}S_1, & \text{otherwise}
				\end{cases}
			\end{align*}
		\end{enumerate}	
	\end{theorem}
	As an application of Theorem \ref{ThIERD}, we set the following result.
	\begin{theorem}
		\begin{enumerate}
			\item [(i)] 
			Risk of the estimator $\delta_{\varphi_1}$ given in (\ref{W}) is nowhere greater than that of $\delta_{{11}}$ under the loss function $L_1(\cdot)$ provided the function $\varphi_1(w)$ satisfies 
			
			\begin{enumerate}
				\item [(a)] $\varphi_1(w)$ is non-decreasing in $w$ and $\lim\limits_{w\rightarrow\infty}\varphi_1(w)=\sqrt{\frac{b(b-1)}{(p_1-1)(p_1-2)}}$
				\item [(b)] $\varphi_1(w)\geq\varphi_*^1(w)$, where $\varphi_*^1(w)$ is $\varphi_{*}^1(w)= \left(\frac{b(b-1) \int_{0}^{\frac{w}{1+w}} s^{p_2-2}(1-s)^{p_1-3} ds} {(p_1+p_2-2)(p_1+p_2-3) \int_{0}^{\frac{w}{1+w}} s^{p_2-2}(1-s)^{p_1-1} ds}\right)^{1/2}$ with $p_1+p_2>3$.
			\end{enumerate}
			
			\item [(ii)] Suppose the following conditions are hold true. 
			\begin{enumerate}
				\item [(a)] $\varphi_1(w)$ is non-decreasing in $w$ and $\lim\limits_{w\rightarrow\infty}\varphi_1(w)=\frac{b-1}{p_1-1}$
				\item [(b)] $\varphi_1(w)\geq\varphi_*^2(w)$, where $\varphi_*^2(w)$ is $\varphi_{*}^2(w) = \frac{ (b-1)\int_{0}^{\frac{w}{1+w}} s^{p_2-2}(1-s)^{p_1-2}ds} {(p_1+p_2-2) \int_{0}^{\frac{w}{1+w}} s^{p_2-2}(1-s)^{p_1-1}ds}$ with $p_1+p_2>2$.
			\end{enumerate}
			Then risk of the estimator $\delta_{\varphi}$ given in (\ref{W}) is nowhere greater than that of $\delta_{12}$ with respect to the loss function $L_2(\cdot)$.	
		\end{enumerate}
	\end{theorem}
	We have derive the improved estimator for $\sigma_2$ under the loss function $L_1(\cdot)$ and $L_2(\cdot)$ respectively.
	
	\begin{theorem}\label{L4.4}
		\begin{enumerate}
			\item[(i)] For $L_1(\cdot)$ loss function, we have  $\psi_{11}(U)=(1+U) \left(\frac{b(b-1)}{(p_1+p_2-2)(p_1+p_2-3)}\right)^{1/2}$, $\psi_{21}(U_1)= (1+p_2U_1) \sqrt{\frac{b(b-1)}{p_2(p_2-1)}}$. The improve estimators of $\sigma_2$ are obtained as follows
			\begin{align*}
				\delta_{21}^1= \max\left\{\psi_{11}(U),\ c_2\right\}S_2~~\mbox{ and }~~
				\delta_{22}^{1} = \begin{cases} 
					\min\left\{\psi_{21}(U_1),\ c_2\right\}S_2, & \ U_1>0\\
					c_2S_2, & \text{otherwise}
				\end{cases}
			\end{align*}
			
			\item [(ii)] Under $L_2(\cdot)$ loss function, we have $\psi_{12}(U)=(1+U)\frac{b-1}{(p_1+p_2-1)}$, $\psi_{22}(U_1)= (1+p_2U_1)\frac{(b-1)}{p_2}$. We get the improve estimators of $\sigma_2$ as 
			\begin{align*}
				& 	\delta_{21}^2 = \max\left\{\psi_{12}(U),\ d_2\right\}S_2~~\mbox{ and } ~~
				\delta_{22}^{2} = \begin{cases} 
					\min\left\{\psi_{22}(U_1),\ d_2\right\}S_2, & \ U_1>0\\
					d_2S_2, & \text{otherwise}
				\end{cases} 
			\end{align*}
			
		\end{enumerate}
	\end{theorem}
	As an application of Theorem \ref{thIERD}, we have the following results.
	\begin{theorem}
		\begin{enumerate}
			\item[(i)] Risk of the estimator  $\delta_{\psi_1}$ given in (\ref{W*}) is nowhere larger than that of $\delta_{12}$ under the loss function $L_1(\cdot)$ provided the function $\psi_1(u)$ satisfies the following conditions
			\begin{enumerate}
				\item [(a)] $\psi_1(u)$ is non-decreasing in $u$ and $\lim\limits_{u\rightarrow 0}\psi_1(u)=\left(\frac{b(b-1)}{(p_2-1)(p_2-2)}\right)^{1/2} $
				\item [(b)] $\psi_1(u)\leq\psi_{*}^1(u)$, where $\psi_{*}^1(u) = \left(\frac{b(b-1) \left[ \Gamma(p_2-2) - \Gamma(p_1+p_2-3) \int_{0}^{\frac{u}{1+u}} s^{p_1-2}(1-s)^{p_2-3} ds \right]} {\Gamma(p_2) - \Gamma(p_1+p_2-1) \int_{0}^{\frac{u}{1+u}} s^{p_1-2}(1-s)^{p_2-1} ds }\right)^{1/2}$ with $p_1+p_2>3$.
			\end{enumerate}
			\item [(ii)] Let us assume that the function $\psi_1(u)$ satisfies the following conditions
			\begin{enumerate}
				\item [(a)] $\psi_1(u)$ is non-decreasing in $u$ and $\lim\limits_{u\rightarrow 0}\psi_1(u)=\frac{b-1}{p_2-1}$
				\item [(b)] $\psi_1(u)\leq\psi_{*}^2(u) $, where $ \psi_{*}^2(u) = \frac{ (b-1) \left[\Gamma(p_2-1) - \Gamma(p_1+p_2-2) \int_{0}^{\frac{u}{1+u}} s^{p_1-2}(1-s)^{p_2-2} ds\right] } { \Gamma(p_2) - \Gamma(p_1+p_2-1) \int_{0}^{\frac{u}{1+u}} s^{p_1-2}(1-s)^{p_2-1} ds}$ provided $p_1+p_2>2$.
			\end{enumerate}			
			Then risk of the estimator $\delta_{\psi_1}$ given in (\ref{W*}) is nowhere larger than that of $\delta_{22}$ under the loss function $L_2(\cdot)$.
		\end{enumerate}
	\end{theorem}
	\subsection{Simulation study}\label{n1}
	In this subsection, we will do a simulation study to compare the risk performance of the estimators proposed in Section \ref{se4.1}. To carry out the simulation study, we have generated 50,000 random samples from two populations having the distributions $Exp\left(\mu_1,\frac{\sigma_1}{\tau}\right)$ and $Exp\left(\mu_2,\frac{\sigma_2}{\tau}\right)$. We have considered various values of $(\mu_1,\mu_2)$ and $(\sigma_1,\sigma_2)$ where $\tau$ from a Gamma distribution $\Gamma(b,1)$ with $b=3,5,7$. Observed that the risk of the estimators depends on the parameters $\sigma_1$ and $\sigma_2$ through $\eta=\sigma_1/\sigma_2$. The performance of the improved estimators is measured by the relative risk improvement (RRI) with respect to the BAEE. The RRI of an estimator $\delta$ is with the respect to $\delta_{1}$ is defined by 
	\begin{equation}\label{rri}
		\mbox{RRI}(\delta)=\frac{\mbox{Risk}(\delta_{1})-\mbox{Risk}(\delta)}{\mbox{Risk}(\delta_{1})}\times 100 \%
	\end{equation}
	In Table \ref{ml1}-\ref{ml4}, we have presented RRI values of $\delta_{12}^2$ and $\delta_{13}^2$ relative to the BAEE $\delta_{21}$ with respect to $L_2(\cdot)$. From the tabulated RRI values, we have observed that the RRI of $\delta_{12}^2$ and $\delta_{13}^2$ are decreasing in $b$ for $b\geq2$. Also the RRI of  $\delta_{12}^2$ and $\delta_{13}^2$ are increasing function of $\eta$ for $0<\eta\leq1$. 
	
	Among the estimators $\delta_{12}^2$ and $\delta_{13}^2$, we observed that $\delta_{13}^2$ outperforms $\delta_{12}^2$ when both $\mu_1$ and $\mu_2$ take negative values (see Table~\ref{ml1}). However, when $\mu_1$ and $\mu_2$ take non-negative values in the neighborhood of zero, $\delta_{12}^2$ performs better than $\delta_{13}^2$ (see Tables~\ref{ml2}--\ref{ml4}).

	In Table \ref{ml5} we have tabulated the RRI $\delta_{11}^2$ and $\delta_{\varphi_*^2}$.
	It is easy to see that the risks of these estimators are independent of the parameters $\mu_1$ and $\mu_2$. Thus we have compute the RRI of $\delta_{11}^2$ and $\delta_{\varphi_*^2}$	only at $(\mu_1,\mu_2)=(0,0)$. Here the RRI of $\delta_{11}^2$ and $\delta_{\varphi_*^2}$ are decreasing function of $b$. RRI of $\delta_{11}^2$ is increasing in $\eta$ but RRI of the estimator $\delta_{\varphi_*^2}$ is not monotone it is increasing when $\eta<0.7$ (approximately) and decreasing when $\eta>0.7$ (approximately). 
	Overall, Table \ref{ml5} shows that  $\delta_{\varphi_*^2}$ performs better than $\delta_{11}^2$ when $0<\eta<0.5$ (approximately) and when $\eta>0.5$ (approximately) then $\delta_{11}^2$ performed better than $\delta_{\varphi_*^2}$.
	
	Now we discuss the numerical performance of the estimators for the parameter $\sigma_2$ under the loss function $L_2(\cdot)$. In Table \ref{ml6}, we present the RRI of $\delta_{21}^2$ and $\delta_{\psi_*^2}$ over $\delta_{22}$ for estimating $\sigma_2$. 
	The RRI of $\delta_{21}$ is increasing function of $\eta$, whereas the RRI of $\delta_{\psi_{*}^2}$ is not monotonic in $\eta$; it is increasing when $\eta<0.7$ (approximately) and then it decreases. 
	Among the estimators $\delta_{21}^2$ and $\delta_{\psi_*^2}$ for estimating $\sigma_2$, we can conclude that  $\delta_{\psi_*^2}$ preferable than $\delta_{21}^2$. A similar observations can be made for the loss function $L_1(\cdot)$.

	\begin{table}[p]
		\centering
		\footnotesize
		\renewcommand{\thetable}{ML1}
		\caption{RRI of improved estimators of $\sigma_1$ under $L_2(\cdot)$ over $\delta_{21}$ when $(\mu_1,\mu_2)=(-0.5,-0.3)$}\label{ml1}
			\begin{tabular}{cc|cc|cc|cc}
				\hline\hline
				$(p_1,p_2)$ & $\eta$
				& \multicolumn{2}{c|}{$b=3$}
				& \multicolumn{2}{c|}{$b=5$}
				& \multicolumn{2}{c}{$b=7$} \\[4pt]
				\cline{3-8} 
				&
				& $\delta_{12}^2$ & $\delta_{13}^2$
				& $\delta_{12}^2$ & $\delta_{13}^2$
				& $\delta_{12}^2$ & $\delta_{13}^2$  \\[4pt]
				\hline
				
				(5,5) & 0.1 & 0.000 & 0.132 & 0.000 & 0.000 & 0.000 & 0.000 \\
				& 0.3 & 0.019 & 0.651 & 0.000 & 0.077 & 0.000 & 0.000 \\
				& 0.5 & 0.055 & 1.004 & 0.037 & 0.146 & 0.000 & 0.002 \\
				& 0.7 & 0.124 & 1.209 & 0.030 & 0.196 & 0.000 & 0.008 \\
				& 0.9 & 0.415 & 1.307& 0.046 & 0.222 & 0.015 & 0.011 \\
				\hline
				
				(5,8) & 0.1 & 0.000 & 0.027 & 0.000 & 0.000 & 0.000 & 0.000\\
				& 0.3 & 0.025 & 0.229 & 0.000 & 0.009 & 0.000 & 0.000 \\
				& 0.5 & 0.061 & 0.397 & 0.045 & 0.027 & 0.000 & 0.002 \\
				& 0.7 & 0.228 & 0.487& 0.041 & 0.042 & 0.000 & 0.002\\
				& 0.9 & 0.514 & 0.528 & 0.068 & 0.047& 0.000 & 0.002\\
				\hline
				
				(15,10) & 0.1 & 0.000 & 0.022  & 0.000 & 0.000 & 0.000 & 0.000 \\
				& 0.3 & 0.000 & 0.133 & 0.000 & 0.000  & 0.000 & 0.000 \\
				& 0.5 & 0.010 & 0.216 & 0.000 & 0.002 & 0.000 & 0.000 \\
				& 0.7 & 0.028 & 0.270 & 0.000 & 0.003 & 0.000 & 0.000\\
				& 0.9 & 0.083 & 0.299 & 0.000 & 0.003 & 0.000 & 0.000 \\
				\hline
				
				(15,21) & 0.1 & 0.000 & 0.000 & 0.000 & 0.000 & 0.000 & 0.000 \\
				& 0.3 & 0.000 & 0.015  & 0.000 & 0.000 & 0.000 & 0.000 \\
				& 0.5 & 0.015 & 0.034& 0.000 & 0.000 & 0.000 & 0.000 \\
				& 0.7 & 0.042 & 0.050 & 0.000 & 0.000 & 0.000 & 0.000 \\
				& 0.9 & 0.123 & 0.058 & 0.000 & 0.000 & 0.000 & 0.000 \\
				\hline\hline
			\end{tabular}
	\end{table}
	%
	\begin{table}[p]
		\centering
		\footnotesize
		\renewcommand{\thetable}{ML2}
		\caption{RRI of improved estimators of $\sigma_1$ under $L_2(\cdot)$ over $\delta_{21}$ when $(\mu_1,\mu_2)=(0,0)$}\label{ml2}
			\begin{tabular}{cc|cc|cc|cc}
				\hline\hline
				$(p_1,p_2)$ & $\eta$
				& \multicolumn{2}{c|}{$b=3$}
				& \multicolumn{2}{c|}{$b=5$}
				& \multicolumn{2}{c}{$b=7$} \\[4pt]
				\cline{3-8} 
				&
				& $\delta_{12}^2$ & $\delta_{13}^2$
				& $\delta_{12}^2$ & $\delta_{13}^2$
				& $\delta_{12}^2$ & $\delta_{13}^2$  \\[4pt]
				\hline
				
				(5,5) & 0.1 & 6.201 & 4.130& 2.488 & 1.299 & 1.395 & 0.574 \\
				& 0.3 & 16.604 & 14.444 &12.089 & 9.706 & 9.943 & 7.432 \\
				& 0.5 & 21.010 & 20.042 & 18.212 & 16.824 & 16.635 & 14.999 \\
				& 0.7 & 22.554 & 22.370 & 21.029 & 20.642 &20.046 & 19.474\\
				& 0.9 & 22.684 & 22.722 &21.806 & 21.869 & 21.168 & 21.172 \\
				\hline
				
				(5,8) & 0.1 & 5.821 & 4.461 & 1.586 & 1.007 & 0.675 & 0.357 \\
				& 0.3 & 18.478 & 16.795 & 12.030 & 10.260 & 9.400 & 7.578 \\
				& 0.5 & 24.330 & 23.598 & 20.102 & 19.043 & 18.303 & 17.029 \\
				& 0.7 & 26.264 & 26.177 & 23.857 & 23.618 &  23.117 & 22.743 \\
				& 0.9 & 26.125 & 26.166 & 24.639 & 24.712 & 24.550 & 24.599 \\
				\hline
				
				(15,10) & 0.1 & 3.511 & 2.845 & 0.538 & 0.311 & 0.108 & 0.055 \\
				& 0.3 & 12.958 & 11.814 & 7.634 & 6.467  & 5.108 & 4.071\\
				& 0.5 & 18.313 & 17.689 & 14.808 & 13.911 & 12.639 & 11.623 \\
				& 0.7 & 20.727 & 20.562& 19.126 & 18.820 & 18.081 & 17.689 \\
				& 0.9 & 21.426 & 21.430 & 20.933 & 20.954 & 20.658 & 20.681 \\
				\hline
				
				(15,21) & 0.1 & 4.173 & 3.737 & 0.375 & 0.266 & 0.073 & 0.058 \\
				& 0.3 & 16.837 & 16.076 & 8.762 & 8.036 & 5.393 & 4.769 \\
				& 0.5 & 24.350 & 23.990 & 18.957 & 18.362 & 15.915 & 15.228 \\
				& 0.7 & 27.291 & 27.280 & 24.972 & 24.840 &  23.813 & 23.565 \\
				& 0.9 & 27.337 & 27.388 & 26.582 & 26.650 & 26.590 & 26.632 \\
				\hline\hline
			\end{tabular}
	\end{table}
	%
	\begin{table}[p]
		\centering
		\footnotesize
		\renewcommand{\thetable}{ML3}
		\caption{RRI of improved estimators of $\sigma_1$ under $L_2(\cdot)$ over $\delta_{21}$ when $(\mu_1,\mu_2)=(0.2,0.5)$}\label{ml3}
			\begin{tabular}{cc|cc|cc|cc}
				\hline\hline
				$(p_1,p_2)$ & $\eta$
				& \multicolumn{2}{c|}{$b=3$}
				& \multicolumn{2}{c|}{$b=5$}
				& \multicolumn{2}{c}{$b=7$} \\[4pt]
				\cline{3-8} 
				&
				& $\delta_{12}^2$ & $\delta_{13}^2$
				& $\delta_{12}^2$ & $\delta_{13}^2$
				& $\delta_{12}^2$ & $\delta_{13}^2$  \\[4pt]
				\hline
				
				(5,5) & 0.1 & 1.594 & 0.354 & 0.0089 & 0.0000 & 0.000 & 0.000 \\
				& 0.3 & 9.077 & 3.727 & 1.7531 & 0.1742 & 0.212 & 0.000 \\
				& 0.5 & 15.690 & 8.247 &5.6334 & 1.0682 & 1.529 & 0.086\\
				& 0.7 & 20.396 & 12.781& 10.4379 & 2.7107 &4.269 & 0.347 \\
				& 0.9 & 23.474 & 16.791 & 15.0131 & 4.9811 & 7.899 & 0.978 \\
				\hline
				
				(5,8) & 0.1 & 2.575 & 0.688 &0.0665 & 0.0000 & 0.000 & 0.000\\
				& 0.3 & 12.724 & 4.553 & 3.2661 & 0.1417 & 0.751 & 0.000\\
				& 0.5 & 20.824 & 9.922 &9.7035 & 1.0646 & 4.095 & 0.090 \\
				& 0.7 & 25.874 & 15.234 & 16.4680 & 2.9339 &9.540 & 0.423 \\
				& 0.9 & 28.599 & 19.912 & 22.0138 & 5.5910 & 15.443 & 1.157\\
				\hline
				
				(15,10) & 0.1 & 0.837 & 0.514 & 0.379 & 0.0000 & 0.000 & 0.000 \\
				& 0.3 & 4.840 & 2.905 & 0.1547 & 0.0197 & 0.000 & 0.000 \\
				& 0.5 & 9.790 & 6.365 & 1.3615 & 0.2858 & 0.075 & 0.000 \\
				& 0.7 & 14.329 & 10.154 & 3.6888 & 1.1976 & 0.516 & 0.050\\
				& 0.9 & 18.041 & 13.771 & 6.8062 & 2.6978 & 1.575 & 0.236 \\
				\hline
				
				(15,21) & 0.1 & 1.849 & 0.544 & 0.0000 & 0.0000 & 0.000 & 0.000 \\
				& 0.3 & 10.270 & 3.983 & 1.2532 & 0.0087 & 0.074 & 0.000 \\
				& 0.5 & 18.835 & 8.858 & 5.9416 & 0.3289 & 1.353 & 0.000\\
				& 0.7 & 25.142 & 14.104 & 12.7033 & 1.5524& 4.866 & 0.059 \\
				& 0.9 & 29.063 & 19.130 & 19.6184 & 3.6281 & 10.365 & 0.304 \\
				\hline\hline
			\end{tabular}
	\end{table}
	%
	%
	\begin{table}[p]
		\centering
		\footnotesize
		\renewcommand{\thetable}{ML4}
		\caption{RRI of improved estimators of $\sigma_1$ under $L_2(\cdot)$ over $\delta_{21}$ when $(\mu_1,\mu_2)=(0.5,1)$}\label{ml4}
			\begin{tabular}{cc|cc|cc|cc}
				\hline\hline
				$(p_1,p_2)$ & $\eta$
				& \multicolumn{2}{c|}{$b=3$}
				& \multicolumn{2}{c|}{$b=5$}
				& \multicolumn{2}{c}{$b=7$} \\[4pt]
				\cline{3-8} 
				&
				& $\delta_{12}^2$ & $\delta_{13}^2$
				& $\delta_{12}^2$ & $\delta_{13}^2$
				& $\delta_{12}^2$ & $\delta_{13}^2$  \\[4pt]
				\hline
				(5,5) & 0.1 & 0.426 & 0.069 & 0.000 & 0.000 & 0.000 & 0.000 \\ 
				~ & 0.3 & 4.155 & 1.454 & 0.189 & 0.000 & 0.000 & 0.000\\ 
				~ & 0.5 & 8.803 & 3.911& 1.143 & 0.142 & 0.085 & 0.000 \\ 
				~ & 0.7 & 13.163 & 6.805 & 2.813 & 0.513 & 0.348 & 0.010\\ 
				~ & 0.9 & 16.928 & 9.747& 5.020 & 1.190 & 0.997 & 0.081 \\ \hline
				(5,8) & 0.1 & 1.213 & 0.282 & 0.000 & 0.000 & 0.000 & 0.000 \\
				~ & 0.3 & 7.493 & 2.039 &0.634 & 0.000 & 0.039 & 0.000  \\ 
				~ & 0.5 & 14.430 & 4.873& 3.013 & 0.115 & 0.489 & 0.000 \\ 
				~ & 0.7 & 20.269 & 8.303 & 6.722 & 0.493 &1.746 & 0.006 \\ 
				~ & 0.9 & 24.689 & 11.791 & 11.026 & 1.234 & 3.941 & 0.086\\ \hline
				(15,10) & 0.1 & 0.335 & 0.217  & 0.000 & 0.000& 0.000 & 0.000 \\
				~ & 0.3 & 1.945 & 1.349 & 0.000 & 0.000 & 0.000 & 0.000\\ 
				~ & 0.5 & 4.359 & 3.086 & 0.059 & 0.015 & 0.000 & 0.000 \\ 
				~ & 0.7 & 7.198 & 5.232 & 0.324 & 0.103 & 0.000 & 0.000\\ 
				~ & 0.9 & 10.085 & 7.616 & 0.960 & 0.340 & 0.018 & 0.000 \\ \hline
				(15,21) & 0.1 & 0.848 & 0.140 & 0.000 & 0.000 & 0.000 & 0.000 \\
				~ & 0.3 & 5.601 & 1.794 & 0.075 & 0.000  & 0.000 & 0.000 \\ 
				~ & 0.5 & 11.790 & 4.245& 1.083 & 0.007 & 0.031 & 0.000 \\ 
				~ & 0.7 & 17.878 & 7.296 & 3.484 & 0.098  & 0.321 & 0.000 \\ 
				~ & 0.9 & 23.063 & 10.626 & 6.941 & 0.401 & 1.232 & 0.000 \\ \hline\hline
			\end{tabular}
	\end{table}
	
	%
	%
	%
	\begin{table}[p]
		\centering
		\footnotesize
		\renewcommand{\thetable}{ML5}
		\caption{RRI values of $\delta_{11}^2$ and $\delta_{\varphi_*^2}$ of $\sigma_1$ under $L_2(\cdot)$ }\label{ml5}

		\begin{tabular}{cc|cc|cc|cc}
			\hline\hline
			$(p_1,p_2)$ & $\eta$
			& \multicolumn{2}{c|}{$b=3$}
			& \multicolumn{2}{c|}{$b=5$}
			& \multicolumn{2}{c}{$b=7$} \\[4pt]
			\cline{3-8} 
			&
			& $\delta_{11}^2$ & $\delta_{\varphi_*^2}$
			& $\delta_{11}^2$ & $\delta_{\varphi_*^2}$
			& $\delta_{11}^2$ & $\delta_{\varphi_*^2}$ \\[4pt]
			\hline
			
			(5,5) & 0.1 & 4.599 & 10.570 & 1.711 & 7.384 & 0.908 & 6.133 \\
			& 0.3 & 13.895 & 17.416 & 9.671 & 15.787 & 7.698 & 14.931 \\
			& 0.5 & 18.277 & 17.207 & 15.538 & 16.203 & 14.023 & 15.509 \\
			& 0.7 & 19.812 & 15.114 & 18.404 & 13.851 & 17.462 & 12.879 \\
			& 0.9 & 19.783 & 12.483 & 19.130 & 10.552 & 18.589 & 9.087 \\
			\hline
			
			(5,8) & 0.1 & 4.849 & 11.486 & 1.205 & 6.947 & 0.477 & 5.494 \\
			& 0.3 & 16.792 & 21.422 & 10.538 & 18.475 & 7.988 & 17.567 \\
			& 0.5 & 22.719 & 21.758 & 18.510 & 19.858 & 16.665 & 19.514 \\
			& 0.7 & 24.633 & 19.026 & 22.351 & 16.792 & 21.600 & 16.266 \\
			& 0.9 & 24.228 & 15.243 & 22.983 & 11.920 & 22.958 & 10.760 \\
			\hline
			
			(15,10) & 0.1 & 2.942 & 5.488 & 0.382 & 1.979 & 0.078 & 1.029 \\
			& 0.3 & 11.668 & 14.917 & 6.596 & 11.596 & 4.284 & 9.971 \\
			& 0.5 & 16.915 & 18.069 & 13.422 & 16.631 & 11.327 & 15.927 \\
			& 0.7 & 19.255 & 18.262 & 17.694 & 17.532 & 16.689 & 17.176 \\
			& 0.9 & 19.792 & 17.132 & 19.400 & 16.165 & 19.170 & 15.517 \\
			\hline
			
			(15,21) & 0.1 & 3.870 & 6.430 & 0.314 & 1.599 & 0.061 & 0.644 \\
			& 0.3 & 16.168 & 19.918 & 8.210 & 14.091 & 4.953 & 11.521 \\
			& 0.5 & 23.664 & 24.917 & 18.232 & 22.277 & 15.218 & 21.193 \\
			& 0.7 & 26.522 & 24.960 & 24.257 & 23.734 & 23.095 & 23.444 \\
			& 0.9 & 26.351 & 22.534 & 25.719 & 20.815 & 25.749 & 20.012 \\
			\hline\hline
			
		\end{tabular}
	\end{table}
	%
	%
	%
	\begin{table}[p]
		\centering
		\footnotesize
		\renewcommand{\thetable}{ML6}
		\caption{RRI of improved estimators of $\sigma_2$ under $L_2(\cdot)$ over $\delta_{22}$ }\label{ml6}
		\begin{tabular}{cc|cc|cc|cc}
			\hline\hline
			$(p_1,p_2)$ & $\eta$
			& \multicolumn{2}{c|}{$b=3$}
			& \multicolumn{2}{c|}{$b=5$}
			& \multicolumn{2}{c}{$b=7$} \\[4pt]
			\cline{3-8} 
			&
			&
			$\delta_{21}^2$ & $\delta_{\psi_*^2}$
			& $\delta_{21}^2$ & $\delta_{\psi_*^2}$
			& $\delta_{21}^2$ & $\delta_{\psi_*^2}$ \\[4pt]
			\hline		
			(5,5) & 0.1 & 0.945 & 8.226 & 0.507 & 8.075 & 0.336 & 8.169 \\
			& 0.3 & 6.764 & 21.075 & 5.657 & 22.857 & 5.315 & 24.027 \\
			& 0.5 & 12.924 & 24.087 & 13.154 & 25.687 & 13.425 & 26.708 \\
			& 0.7 & 17.546 & 20.372 & 19.624 & 19.528 & 20.896 & 19.100 \\
			& 0.9 & 20.203 & 11.894 & 23.876 & 6.790 & 26.072 & 3.858 \\
			\hline
			
			(5,8) & 0.1 & 0.397 & 4.467 & 0.116 & 4.115 & 0.054 & 4.029 \\
			& 0.3 & 4.067 & 14.341 & 2.911 & 15.182 & 2.439 & 15.616 \\
			& 0.5 & 8.816 & 19.149 & 8.274 & 20.364 & 7.929 & 20.676 \\
			& 0.7 & 12.836 & 19.738 & 13.666 & 19.851 & 13.983 & 19.062 \\
			& 0.9 & 15.612 & 17.117 & 17.802 & 14.711 & 18.817 & 11.919 \\
			\hline
			
			(15,10) & 0.1 & 0.508 & 2.247 & 0.114 & 1.316 & 0.025 & 0.892 \\
			& 0.3 & 5.660 & 12.861 & 4.100 & 13.232 & 3.048 & 13.041 \\
			& 0.5 & 11.867 & 19.625 & 12.022 & 23.145 & 11.312 & 24.602 \\
			& 0.7 & 15.991 & 20.208 & 19.143 & 24.702 & 20.098 & 26.399 \\
			& 0.9 & 17.086 & 15.269 & 22.600 & 17.661 & 25.090 & 17.424 \\
			\hline
			
			(15,21) & 0.1 & 0.293 & 1.012 & 0.037 & 0.413 & 0.005 & 0.214 \\
			& 0.3 & 3.708 & 7.811 & 2.207 & 6.970 & 1.398 & 6.410 \\
			& 0.5 & 8.723 & 14.543 & 7.774 & 15.459 & 6.824 & 15.914 \\
			& 0.7 & 13.056 & 18.768 & 14.027 & 20.645 & 14.005 & 21.798 \\
			& 0.9 & 15.732 & 20.344 & 18.542 & 21.150 & 19.689 & 21.701 \\
			\hline\hline
		\end{tabular}
	\end{table}

			\subsection{Exponential inverse Gaussian (E-IG) distribution} \label{se4.2}
			In this subsection, we will derive the improved estimator of the parameter for the exponential inverse Gaussian distribution under the loss functions $L_1(\cdot)$ and $L_2(\cdot)$. In this case, we have $E(\tau)=m$, $E(\tau^{-1})=\frac{1}{m}+\frac{1}{n}$ (see \cite{seshadri1994inverse}). We get the BAEE of $\sigma_i$ under the loss function $L_1(\cdot)$ is as $\delta_{1i}=c_{i}S_i$ with $c_{i}=\left(\frac{m}{(p_i-1)(p_i-2)\left(\frac{1}{m}+\frac{1}{n}\right)}\right)^{1/2}$ for $i=1,2$ and for the loss function $L_2(\cdot)$ we have the BAEE of $\sigma_i$ is $\delta_{{2i}}=d_iS_i$ with $d_{i}=\frac{1}{(p_i-1)\left(\frac{1}{m}+\frac{1}{n}\right)}.$
			As an application of Theorem \ref{Th1}, \ref{ths21} and \ref{th2.4}, we obtained the following results as follows.      

			\begin{theorem}
				\begin{enumerate}
					\item [(i)] Under the loss function $L_1(\cdot)$, we have $\varphi_{11}(W)= \frac{m^{1/2}(1+W)}{\left((p_1+p_2-2)(p_1+p_2-3)\left(1/m+1/n\right)\right)^{1/2}}$, $\varphi_{21}(W,W_1)=\frac{(1+W+p_1W_1) m^{1/2} }{\left((p_1+p_2-1)(p_1+p_2-2)\left(1/m+1/n\right)\right)^{1/2}}$, $\varphi_{31}(W,W_2)=\frac{(1+W+p_2W_2) m^{1/2} }{\left((p_1+p_2-1)(p_1+p_2-2)\left(1/m+1/n\right)\right)^{1/2}}$. The improved estimator of $\sigma_1$ are obtained as follows
					\begin{align*}
						&  \delta^1_{11}=\min\left\{\varphi_{11}(W),c_1\right\}S_1,~~~
						\delta_{12}^{1} = \begin{cases} 
							\min\left\{\varphi_{21}(W,W_1),c_{1}\right\}S_1, & \ W_1>0\\
							c_{1}S_1, & \text{otherwise}
						\end{cases}\\
						&  	\delta_{13}^{1} = \begin{cases} 
							\min\left\{\varphi_{31}(W,W_2),c_{1}\right\}S_1, & \ W_2>0\\
							c_{1}S_1, & \text{otherwise}
						\end{cases}
					\end{align*}
					\item [(ii)] For the loss function $L_2(\cdot)$, we have $\varphi_{21}(W)=\frac{(1+W)}{(p_1+p_2-2)\left(\frac{1}{m}+\frac{1}{n}\right)}$, $\varphi_{22}(W,W_1)=\frac{(1+W+p_1W_1)}{(p_1+p_2-1)\left(\frac{1}{m}+\frac{1}{n}\right)}$, $\varphi_{32}(W,W_1)=\frac{(1+W+p_2W_2)}{(p_1+p_2-1)(1/m+1/n)}$. We have the improved estimator of $\sigma_1$ are as follows
					\begin{align*}
						& \delta^2_{11} = \min\left\{\varphi_{21}(W),d_1\right\}S_1,~~
						\delta_{12}^{2} = \begin{cases} 
							\min\left\{\varphi_{22}(W,W_1),d_{1}\right\}S_1, &\ W_1>0\\
							d_{1}S_1, & \text{otherwise}
						\end{cases}\\
						& 	\delta_{13}^{2} = \begin{cases} 
							\min\left\{\varphi_{32}(W,W_2),d_{1}\right\}S_1, & \ W_2>0\\
							d_{1}S_1, & \text{otherwise}
						\end{cases}
					\end{align*}
				\end{enumerate}
			\end{theorem}
			Applying Theorem \ref{thIERD}, we obtained the results as follows.
			\begin{theorem}
				\begin{enumerate}
					\item[(i)] Under the loss function $L_1(\cdot)$, the risk of the  estimator $\delta_{\varphi_1}$ is nowhere larger than that of $\delta_{11}$ provided the function $\varphi_1(w)$ satisfies
					\begin{enumerate}
						\item [(a)] $\varphi_1(w)$ is non-decreasing in $w$ and $\lim\limits_{w\rightarrow\infty}\varphi_1(w)=\left(\frac{m}{(p_1-1)(p_1-2)\left(1/m+ 1/n \right)}\right)^{1/2}$
						\item [(b)] $\varphi_1(w)\geq \varphi_{*}^1(w)= \left(\frac{m \int_{0}^{\frac{w}{1+w}} s^{p_2-2}(1-s)^{p_1-3} ds} {\left(1/m+ 1/n\right)(p_1+p_2-2)(p_1+p_2-3) \int_{0}^{\frac{w}{1+w}} s^{p_2-2}(1-s)^{p_1-1} ds}\right)^{1/2}$ with $p_1+p_2>3$.
					\end{enumerate}
					\item [(ii)] Let us assume that the function  $\varphi_1(w)$ satisfies the following conditions:
					\begin{enumerate}
						\item [(a)] $\varphi_1(w)$ is non-decreasing in $w$ and $\lim\limits_{w\rightarrow\infty}\varphi_1(w)=\frac{1}{(p_1-1)\left(1/m + 1/n\right)}$
						\item [(b)] $\varphi_1(w)\geq \varphi_{*}^2(w) = \frac{ \int_{0}^{\frac{w}{1+w}} s^{p_2-2}(1-s)^{p_1-2}ds} {\left(1/m + 1/n\right)(p_1+p_2-2) \int_{0}^{\frac{w}{1+w}} s^{p_2-2}(1-s)^{p_1-1}ds}$ with $p_1+p_2>2$	
						
						
					\end{enumerate}
					Then risk of the estimator $\delta_{\varphi_1}$ is nowhere larger than that of $\delta_{21}$ under the loss function $L_2(\cdot)$.
				\end{enumerate}
			\end{theorem}
			
			Now we have obtained the improve estimator of $\sigma_2$ with respect to the loss functions $L_1(\cdot)$ and  $L_2(\cdot)$ respectively.
			
			\begin{theorem}
				\begin{enumerate}
					\item[(i)] For the loss $L_1(\cdot)$, we have $\psi_{11}(U)=(1+U) \left(\frac{m}{(p_1+p_2-2)(p_1+p_2-3)(1/m+1/n)}\right)^{1/2}$, $\psi_{21}(U_1)= \frac{(1+p_2U_1)m^{1/2} }{\left(p_2(p_2-1)(1/m+1/n)\right)^{1/2}}$. The improve estimator of $\sigma_2$ are obtained as follows
					\begin{align*}
						\delta_{21}^1 = \max\left\{\psi_{11}(U),c_2\right\}S_2 ~~\mbox{ and }~~
						\delta_{22}^{1} = \begin{cases} 
							\min\left\{\psi_{21}(U_1),c_2\right\}S_2, & \ U_1>0\\
							c_2S_2, & \text{otherwise}
						\end{cases}
					\end{align*}
					\item [(ii)] Under $L_2(\cdot)$ loss function, we have $\psi_{12}(U)=\frac{(1+U)}{(p_1+p_2-2)(1/m+1/n)}$, $\psi_{22}(U_1)= \frac{(1+p_2U_1)}{p_2(1/m+1/n)}$. We get the improve estimator of $\sigma_2$ as
					\begin{align*}
						&	\delta_{21}^2 = \max\left\{\psi_{12}(U),d_2\right\}S_2~~\mbox{ and }~~
						\delta_{22}^{2} = \begin{cases} 
							\min\left\{\psi_{22}(U_1),d_2\right\}S_2, & \ U_1>0\\
							d_2S_2, & \text{otherwise}
						\end{cases}
					\end{align*}
					
				\end{enumerate}
			\end{theorem} 
			Applying Theorem \ref{ThIERD}, we obtained the results as follows.
			\begin{theorem}
				\begin{enumerate}
					\item [(i)] Risk of the estimator $\delta_{\psi_1}$ is nowhere larger than that of $\delta_{12}$ under the loss function $L_1(\cdot)$ provided the function $\psi_1(u)$ satisfies the following conditions
					\begin{enumerate}
						\item [(a)] $\psi_1(u)$ is non-decreasing in $u$ and $\lim\limits_{u\rightarrow 0}\psi_1(u)=\left(\frac{m}{(p_2-1)(p_2-2)(1/m+1/n)}\right)^{1/2} $
						\item [(b)] $\psi_1(u)\leq \psi_{*}^1(u) = \left(\frac{m \left[ \Gamma(p_2-2) - \Gamma(p_1+p_2-3) \int_{0}^{\frac{u}{1+u}} s^{p_1-2}(1-s)^{p_2-3} ds \right]} {\left(1/m + 1/n\right) \left[\Gamma(p_2) - \Gamma(p_1+p_2-1) \int_{0}^{\frac{u}{1+u}} s^{p_1-2}(1-s)^{p_2-1} ds \right]}\right)^{1/2}$ provided $p_1+p_2>3$.

						
					\end{enumerate}
					\item [(ii)] Let us assume that the function $\psi_{1}(u)$ satisfies the following conditions
					\begin{enumerate}
						\item [(a)] $\psi_1(u)$ is non-decreasing in $u$ and $\lim\limits_{u\rightarrow 0}\psi_1(u)=\frac{1}{(p_2-1)(1/m+1/n)}$
						\item [(b)] $\psi_1(u)\leq  \psi_{*}^2(u) = \frac{  \Gamma(p_2-1) - \Gamma(p_1+p_2-2) \int_{0}^{\frac{u}{1+u}} s^{p_1-2}(1-s)^{p_2-2} ds } {\left(1/m + 1/n\right) \left[\Gamma(p_2) - \Gamma(p_1+p_2-1) \int_{0}^{\frac{u}{1+u}} s^{p_1-2}(1-s)^{p_2-1} ds \right]}$ provided $p_1+p_2>2$. 
						
						
					\end{enumerate}		
					Then risk of the estimator $\delta_{\psi_1}$ defined in (\ref{W*}) is nowhere larger than that of $\delta_{22}$.
				\end{enumerate}
			\end{theorem}
			
			\subsection{Simulation study}
			In this subsection, we have discussed the risk performance of the estimators derived in Section \ref{se4.2} using simulations. For this purpose, we generate the mixing distribution from an inverse Gaussian distribution $IG(m,n)$ with parameter values $(m,n)=(3,5),(4,6),(8,10)$. The performance of the improved estimators is measured using the RRI with respect to the BAEE (see \ref{rri}). We discuss the risk performance of $\sigma_1$ and $\sigma_2$ for the E-IG distribution under the loss function $L_2(\cdot)$. In Table \ref{g1}-\ref{g5}, we have presented the RRI values of $\delta_{12}^2$ and $\delta_{13}^2$ with the respect to the BAEE $\delta_{21}$ for the loss function $L_2(\cdot)$. From the tabulated RRi values, we observe that the RRI of $\delta_{12}^2$ and $\delta_{13}^2$ is an increasing function of $\eta$ for $0<\eta\leq1$. Also the RRI of $\delta_{12}^2$ and $\delta_{13}^2$ decrease as $m$ and $n$ increases.
			
			Among the estimators $\delta_{12}^2$ and $\delta_{13}^2$, we observed that $\delta_{13}^2$ outperforms $\delta_{12}^2$ when both $\mu_1$ and $\mu_2$ take negative values (see Table~\ref{g1}). However, when $\mu_1$ and $\mu_2$ take non-negative values in the neighborhood of zero, $\delta_{12}^2$ performs better than $\delta_{13}^2$ (see Tables~\ref{g2}--\ref{g4}).
			
			In Table \ref{g5}, we have tabulated the RRI of $\delta_{11}^2$ and $\delta_{\varphi_*^2}$. It is easy to see that the risks of these estimators are independent of the parameters $\mu_1$ and $\mu_2$. So without loss of generality we take $(\mu_1,\mu_2)=(0,0)$. Here the $\delta_{11}^2$ and $\delta_{\varphi_*^2}$ are decreasing function of $m$ and $n$. RRI of $\delta_{11}^2$ is an increasing function $\eta$ but RRI of $\delta_{\varphi_*^2}$ is not monotone function of $\eta$, it is increasing when $\eta < 0.6$ (approximately) and decreases for $\eta > 0.6$ (approximately). Overall, Table \ref{g5} shows that $\delta_{\varphi_*^2}$ performs better than $\delta_{11}^2$ when $0<\eta<0.6$ (approximately) and when $\eta>0.6$ (approximately) then $\delta_{11}^2$ performs better than $\delta_{\varphi_*^2}$. Now we discuss the risk performance of the estimators of $\sigma_2$ under the loss function $L_2(\cdot)$. In Table \ref{g6}, we present the RRI of $\delta_{21}^2$ and $\delta_{\psi_*^2}$ over $\delta_{22}$ for estimating $\sigma_2$. The RRI of $\delta_{21}^2$ is increasing function of $\eta$ whereas the RRI of $\delta_{\psi_*^2}$ is not monotonic in $\eta$; it increases for $0<\eta<0.7$ (approximately) and decreases for $\eta>0.7$ (approximately). The RRI of both estimators $\delta_{21}^2$ and $\delta_{\psi_*^2}$ are increasing function of $(m,n)$. In Table \ref{g6}, the estimator $\delta_{\psi_{*}^2}$ performs better than $\delta_{21}^2$. 
			Among the estimators $\delta_{21}^2$ and $\delta_{\psi_*^2}$ for estimating $\sigma_2$, we can observe that $\delta_{\psi_*^2}$ is preferable than $\delta_{21}^2$. The same observations can be made for the loss function $L_1(\cdot)$.

			\begin{table}[p]
				\centering
				\footnotesize
				\renewcommand{\thetable}{G1}
				\caption{RRI of improved estimators of $\sigma_1$ under $L_2(\cdot)$ over $\delta_{21}$ when $(\mu_1,\mu_2)=(-0.5,-0.3)$}\label{g1}
					\begin{tabular}{cc|cc|cc|cc}
						\hline\hline
						& $(m,n)$ 
						& \multicolumn{2}{c|}{(3,5)}
						& \multicolumn{2}{c|}{(4,6)}
						& \multicolumn{2}{c}{(8,10)} \\[4pt]
						\hline
						$(p_1,p_2)$ & $\eta$
						& $\delta_{12}^2$ & $\delta_{13}^2$
						& $\delta_{12}^2$ & $\delta_{13}^2$
						& $\delta_{12}^2$ & $\delta_{13}^2$ \\[4pt]
						\hline			
						(5,5) & 0.1 & 0.000 & 0.000 & 0.000 & 0.000 & 0.000 & 0.000 \\
						~ & 0.3 & 0.000 & 0.243 & 0.000 & 0.165 & 0.000 & 0.068 \\
						~ & 0.5 & 0.102 & 0.993 & 0.057 & 0.702 & 0.022 & 0.254 \\
						~ & 0.7 & 0.435 & 1.922 & 0.248 & 1.381 & 0.045 & 0.512 \\
						~ & 0.9 & 1.169 & 2.879 & 0.670 & 2.097 & 0.153 & 0.783 \\
						\hline			
						(5,8) & 0.1 & 0.000 & 0.001 & 0.000 & 0.000 & 0.000 & 0.000 \\
						~ & 0.3 & 0.000 & 0.100 & 0.000 & 0.047 & 0.000 & 0.012 \\
						~ & 0.5 & 0.080 & 0.448 & 0.030 & 0.273 & 0.005 & 0.071 \\
						~ & 0.7 & 0.439 & 0.963 & 0.271 & 0.624 & 0.031 & 0.169 \\
						~ & 0.9 & 1.232 & 1.571 & 0.735 & 1.029 & 0.137 & 0.274 \\
						\hline			
						(15,10) & 0.1 & 0.000 & 0.000 & 0.000 & 0.000 & 0.000 & 0.000 \\
						~ & 0.3 & 0.000 & 0.004 & 0.000 & 0.001 & 0.000 & 0.000 \\
						~ & 0.5 & 0.001 & 0.067 & 0.000 & 0.030 & 0.000 & 0.003 \\
						~ & 0.7 & 0.006 & 0.263 & 0.007 & 0.142 & 0.000 & 0.011 \\
						~ & 0.9 & 0.019 & 0.560 & 0.011 & 0.297 & 0.000 & 0.035 \\
						\hline			
						(15,21) & 0.1 & 0.000 & 0.000 & 0.000 & 0.000 & 0.000 & 0.000 \\
						~ & 0.3 & 0.000 & 0.000 & 0.000 & 0.000 & 0.000 & 0.000 \\
						~ & 0.5 & 0.000 & 0.006 & 0.000 & 0.000 & 0.000 & 0.000 \\
						~ & 0.7 & 0.004 & 0.029 & 0.005 & 0.007 & 0.000 & 0.000 \\
						~ & 0.9 & 0.022 & 0.074 & 0.011 & 0.031 & 0.000 & 0.000 \\
						\hline\hline
					\end{tabular}
			\end{table}

			\begin{table}[p]
				\centering
				\footnotesize
				\renewcommand{\thetable}{G2}
				\caption{RRI of improved estimators of $\sigma_1$ under $L_2(\cdot)$ over $\delta_{21}$ when $(\mu_1,\mu_2)=(0,0)$}\label{g2}
				\begin{tabular}{cc|cc|cc|cc}
					\hline\hline
					& $(m,n)$
					& \multicolumn{2}{c|}{(3,5)}
					& \multicolumn{2}{c|}{(4,6)}
					& \multicolumn{2}{c}{(8,10)} \\[4pt]
					\hline
					$(p_1,p_2)$ & $\eta$
					& $\delta_{12}^2$ & $\delta_{13}^2$
					& $\delta_{12}^2$ & $\delta_{13}^2$
					& $\delta_{12}^2$ & $\delta_{13}^2$ \\[4pt]
					\hline			
					(5,5) & 0.1 & 0.097 & 0.028 & 0.093 & 0.026 & 0.083 & 0.023 \\
					~ & 0.3 & 1.670 & 1.043 & 1.575 & 0.978 & 1.433 & 0.897 \\
					~ & 0.5 & 4.055 & 3.363 & 3.826 & 3.163 & 3.503 & 2.914 \\
					~ & 0.7 & 5.808 & 5.433 & 5.469 & 5.108 & 5.011 & 4.700 \\
					~ & 0.9 & 6.632 & 6.466 & 6.235 & 6.069 & 5.707 & 5.580 \\
					\hline			
					(5,8) & 0.1 & 0.012 & 0.002 & 0.011 & 0.002 & 0.010 & 0.002 \\
					~ & 0.3 & 1.068 & 0.739 & 1.006 & 0.696 & 0.913 & 0.635 \\
					~ & 0.5 & 4.028 & 3.434 & 3.807 & 3.236 & 3.483 & 2.966 \\
					~ & 0.7 & 6.815 & 6.474 & 6.429 & 6.106 & 5.893 & 5.614 \\
					~ & 0.9 & 8.252 & 8.122 & 7.774 & 7.644 & 7.120 & 7.022 \\
					\hline			
					(15,10) & 0.1 & 0.000 & 0.000 & 0.000 & 0.000 & 0.000 & 0.000 \\
					~ & 0.3 & 0.075 & 0.050 & 0.070 & 0.047 & 0.061 & 0.041 \\
					~ & 0.5 & 0.601 & 0.484 & 0.558 & 0.451 & 0.489 & 0.398 \\
					~ & 0.7 & 1.543 & 1.426 & 1.439 & 1.337 & 1.267 & 1.187 \\
					~ & 0.9 & 2.228 & 2.151 & 2.080 & 2.012 & 1.827 & 1.783 \\
					\hline			
					(15,21) & 0.1 & 0.000 & 0.000 & 0.000 & 0.000 & 0.000 & 0.000 \\
					~ & 0.3 & 0.010 & 0.006 & 0.010 & 0.006 & 0.008 & 0.005 \\
					~ & 0.5 & 0.358 & 0.315 & 0.331 & 0.290 & 0.287 & 0.252 \\
					~ & 0.7 & 1.725 & 1.637 & 1.612 & 1.532 & 1.412 & 1.348 \\
					~ & 0.9 & 3.304 & 3.269 & 3.090 & 3.060 & 2.719 & 2.702 \\
					\hline\hline
				\end{tabular}	 
			\end{table}
			\begin{table}[p]
				\centering
				\footnotesize
				\centering
				\renewcommand{\thetable}{G3}
				\caption{table}{RRI of improved estimators of $\sigma_1$ under $L_2(\cdot)$ over $\delta_{21}$ when $(\mu_1,\mu_2)=(0.2,0.5)$}\label{g3}
				
				\begin{tabular}{cc|cc|cc|cc}
					\hline\hline
					
					& $(m,n)$
					& \multicolumn{2}{c|}{(3,5)}
					& \multicolumn{2}{c|}{(4,6)}
					& \multicolumn{2}{c}{(8,10)} \\[4pt]
					\hline
					
					$(p_1,p_2)$ & $\eta$
					& $\delta_{12}^2$ & $\delta_{13}^2$
					& $\delta_{12}^2$ & $\delta_{13}^2$
					& $\delta_{12}^2$ & $\delta_{13}^2$ \\[4pt]
					\hline
					
					(5,5) & 0.1 & 0.001 & 0.000 & 0.000 & 0.000 & 0.000 & 0.000 \\
					~ & 0.3 & 0.483 & 0.041 & 0.330 & 0.019 & 0.083 & 0.000 \\
					~ & 0.5 & 2.630 & 0.560 & 2.031 & 0.340 & 0.847 & 0.066 \\
					~ & 0.7 & 5.935 & 2.182 & 4.942 & 1.479 & 2.625 & 0.420 \\
					~ & 0.9 & 9.297 & 4.753 & 8.164 & 3.479 & 5.050 & 1.252 \\
					\hline
					
					(5,8) & 0.1 & 0.000 & 0.000 & 0.000 & 0.000 & 0.000 & 0.000 \\
					~ & 0.3 & 0.439 & 0.017 & 0.333 & 0.006 & 0.145 & 0.000 \\
					~ & 0.5 & 2.951 & 0.493 & 2.445 & 0.314 & 1.310 & 0.061 \\
					~ & 0.7 & 7.138 & 2.113 & 6.304 & 1.435 & 4.092 & 0.427 \\
					~ & 0.9 & 11.290 & 5.012 & 10.458 & 3.648 & 7.730 & 1.327 \\
					\hline
					
					(15,10) & 0.1 & 0.000 & 0.000 & 0.000 & 0.000 & 0.000 & 0.000 \\
					~ & 0.3 & 0.001 & 0.000 & 0.000 & 0.000 & 0.000 & 0.000 \\
					~ & 0.5 & 0.098 & 0.024 & 0.055 & 0.010 & 0.009 & 0.000 \\
					~ & 0.7 & 0.693 & 0.229 & 0.453 & 0.126 & 0.115 & 0.022 \\
					~ & 0.9 & 2.241 & 1.002 & 1.605 & 0.612 & 0.513 & 0.130 \\
					\hline
					
					(15,21) & 0.1 & 0.000 & 0.000 & 0.000 & 0.000 & 0.000 & 0.000 \\
					~ & 0.3 & 0.000 & 0.000 & 0.000 & 0.000 & 0.000 & 0.000 \\
					~ & 0.5 & 0.118 & 0.002 & 0.082 & 0.000 & 0.022 & 0.000 \\
					~ & 0.7 & 1.155 & 0.136 & 0.870 & 0.072 & 0.351 & 0.003 \\
					~ & 0.9 & 4.213 & 0.874 & 3.476 & 0.514 & 1.748 & 0.092 \\
					\hline\hline
				\end{tabular}
			\end{table}

			\begin{table}[p]
				\centering
				\footnotesize
				\renewcommand{\thetable}{G4}
				\caption{RRI of improved estimators of $\sigma_1$ under $L_2(\cdot)$ over $\delta_{21}$ when $(\mu_1,\mu_2)=(0.5,1)$}\label{g4}
				\begin{tabular}{cc|cc|cc|cc}
					\hline\hline
					
					& $(m,n)$
					& \multicolumn{2}{c|}{(3,5)}
					& \multicolumn{2}{c|}{(4,6)}
					& \multicolumn{2}{c}{(8,10)} \\[4pt]
					\hline
					
					$(p_1,p_2)$ & $\eta$
					& $\delta_{12}^2$ & $\delta_{13}^2$
					& $\delta_{12}^2$ & $\delta_{13}^2$
					& $\delta_{12}^2$ & $\delta_{13}^2$ \\[4pt]
					\hline
					
					(5,5) & 0.1 & 0.000 & 0.000 & 0.000 & 0.000 & 0.000 & 0.000 \\
					~ & 0.3 & 0.060 & 0.001 & 0.027 & 0.000 & 0.000 & 0.000 \\
					~ & 0.5 & 0.750 & 0.097 & 0.469 & 0.046 & 0.082 & 0.000 \\
					~ & 0.7 & 2.445 & 0.573 & 1.649 & 0.306 & 0.471 & 0.036 \\
					~ & 0.9 & 4.910 & 1.639 & 3.597 & 0.993 & 1.264 & 0.176 \\
					\hline
					
					(5,8) & 0.1 & 0.000 & 0.000 & 0.000 & 0.000 & 0.000 & 0.000 \\
					~ & 0.3 & 0.127 & 0.000 & 0.073 & 0.000 & 0.011 & 0.000 \\
					~ & 0.5 & 1.255 & 0.091 & 0.878 & 0.039 & 0.297 & 0.002 \\
					~ & 0.7 & 4.070 & 0.580 & 3.080 & 0.317 & 1.238 & 0.042 \\
					~ & 0.9 & 7.906 & 1.720 & 6.368 & 1.062 & 3.057 & 0.201 \\
					\hline
					
					(15,10) & 0.1 & 0.000 & 0.000 & 0.000 & 0.000 & 0.000 & 0.000 \\
					~ & 0.3 & 0.000 & 0.000 & 0.000 & 0.000 & 0.000 & 0.000 \\
					~ & 0.5 & 0.005 & 0.000 & 0.000 & 0.000 & 0.000 & 0.000 \\
					~ & 0.7 & 0.086 & 0.032 & 0.039 & 0.011 & 0.000 & 0.000 \\
					~ & 0.9 & 0.413 & 0.182 & 0.220 & 0.087 & 0.027 & 0.005 \\
					\hline
					
					(15,21) & 0.1 & 0.000 & 0.000 & 0.000 & 0.000 & 0.000 & 0.000 \\
					~ & 0.3 & 0.000 & 0.000 & 0.000 & 0.000 & 0.000 & 0.000 \\
					~ & 0.5 & 0.016 & 0.000 & 0.005 & 0.000 & 0.000 & 0.000 \\
					~ & 0.7 & 0.319 & 0.006 & 0.186 & 0.000 & 0.032 & 0.000 \\
					~ & 0.9 & 1.649 & 0.134 & 1.076 & 0.057 & 0.269 & 0.000 \\
					\hline\hline
				\end{tabular}
			\end{table}
			
			
			\begin{table}[p]
				\centering
				\footnotesize
				\centering
				\renewcommand{\thetable}{G5}
				\caption{table}{RRI values of $\delta_{11}^2$ and $\delta_{\varphi_*^2}$ of $\sigma_1$ under $L_2(\cdot)$}
				\label{g5}
				
				\begin{tabular}{cc|cc|cc|cc}
					\hline\hline
					& $(m,n)$
					& \multicolumn{2}{c|}{(3,5)}
					& \multicolumn{2}{c|}{(4,6)}
					& \multicolumn{2}{c}{(8,10)}\\[4pt]
					\hline
					$(p_1,p_2)$ & $\eta$
					& $\delta_{11}^2$ & $\delta_{\varphi_*^2}$
					& $\delta_{11}^2$ & $\delta_{\varphi_*^2}$
					& $\delta_{11}^2$ & $\delta_{\varphi_*^2}$\\[4pt]
					\hline
					
					(5,5) & 0.1 & 0.054 & 1.304 & 0.052 & 1.228 & 0.046 & 1.125 \\
					~ & 0.3 & 1.209 & 4.437 & 1.137 & 4.166 & 1.039 & 3.821 \\
					~ & 0.5 & 3.217 & 4.616 & 3.026 & 4.311 & 2.778 & 3.950 \\
					~ & 0.7 & 4.854 & 3.163 & 4.560 & 2.915 & 4.193 & 2.667 \\
					~ & 0.9 & 5.651 & 1.028 & 5.296 & 0.879 & 4.865 & 0.797 \\
					\hline
					(5,8) & 0.1 & 0.006 & 0.947 & 0.006 & 0.894 & 0.005 & 0.817 \\
					~ & 0.3 & 0.815 & 5.185 & 0.768 & 4.885 & 0.698 & 4.471 \\
					~ & 0.5 & 3.457 & 6.308 & 3.260 & 5.923 & 2.980 & 5.419 \\
					~ & 0.7 & 6.173 & 4.741 & 5.822 & 4.416 & 5.343 & 4.034 \\
					~ & 0.9 & 7.671 & 1.757 & 7.217 & 1.569 & 6.623 & 1.424 \\
					\hline
					(15,10) & 0.1 & 0.000 & 0.011 & 0.000 & 0.011 & 0.000 & 0.009 \\
					~ & 0.3 & 0.058 & 0.690 & 0.054 & 0.643 & 0.046 & 0.567 \\
					~ & 0.5 & 0.512 & 1.647 & 0.475 & 1.534 & 0.417 & 1.349 \\
					~ & 0.7 & 1.362 & 1.649 & 1.273 & 1.527 & 1.124 & 1.325 \\
					~ & 0.9 & 2.024 & 0.703 & 1.892 & 0.634 & 1.671 & 0.513 \\
					\hline
					(15,21) & 0.1 & 0.000 & 0.001 & 0.000 & 0.001 & 0.000 & 0.001 \\
					~ & 0.3 & 0.008 & 0.456 & 0.007 & 0.425 & 0.006 & 0.373 \\
					~ & 0.5 & 0.320 & 2.096 & 0.295 & 1.955 & 0.256 & 1.722 \\
					~ & 0.7 & 1.627 & 2.837 & 1.522 & 2.640 & 1.334 & 2.316 \\
					~ & 0.9 & 3.211 & 1.554 & 3.005 & 1.433 & 2.651 & 1.224 \\
					\hline\hline
				\end{tabular}
			\end{table}

			\begin{table}[p]
				\centering
				\footnotesize
				\renewcommand{\thetable}{G6}
				\caption{table}{RRI of improved estimators of $\sigma_2$ under $L_2(\cdot)$ over $\delta_{22}$}
				\label{g6}
				
				\begin{tabular}{cc|cc|cc|cc}
					\hline\hline
					& $(m,n)$
					& \multicolumn{2}{c|}{(3,5)}
					& \multicolumn{2}{c|}{(4,6)}
					& \multicolumn{2}{c}{(8,10)}\\[4pt]
					\hline
					$(p_1,p_2)$ & $\eta$
					& $\delta_{21}^2$ & $\delta_{\psi_*^2}$
					& $\delta_{21}^2$ & $\delta_{\psi_*^2}$
					& $\delta_{21}^2$ & $\delta_{\psi_*^2}$\\[4pt]
					\hline
					(5,5) & 0.1 & 2.050 & 10.377 & 2.292 & 10.642 & 2.726 & 11.113 \\
					& 0.3 & 12.451 & 25.715 & 12.926 & 25.882 & 13.775 & 26.209 \\
					& 0.5 & 21.428 & 29.647 & 21.781 & 29.855 & 22.435 & 30.276 \\
					& 0.7 & 27.048 & 26.216 & 27.230 & 26.693 & 27.627 & 27.588 \\
					& 0.9 & 29.401 & 17.703 & 29.478 & 18.643 & 29.725 & 20.330 \\
					\hline
					(5,8) & 0.1 & 1.020 & 5.968 & 1.201 & 6.198 & 1.546 & 6.639 \\
					& 0.3 & 8.523 & 18.297 & 8.994 & 18.545 & 9.862 & 19.064 \\
					& 0.5 & 16.189 & 24.404 & 16.628 & 24.670 & 17.485 & 25.239 \\
					& 0.7 & 21.829 & 25.709 & 22.160 & 26.123 & 22.868 & 26.945 \\
					& 0.9 & 25.161 & 23.428 & 25.419 & 24.125 & 26.029 & 25.411 \\
					\hline
					(15,10) & 0.1 & 1.111 & 3.762 & 1.337 & 4.124 & 1.833 & 4.845 \\
					& 0.3 & 11.309 & 19.323 & 11.976 & 19.751 & 13.267 & 20.626 \\
					& 0.5 & 21.588 & 28.613 & 22.100 & 28.777 & 23.131 & 29.238 \\
					& 0.7 & 27.706 & 30.255 & 27.940 & 30.368 & 28.538 & 30.749 \\
					& 0.9 & 29.194 & 25.554 & 29.272 & 25.918 & 29.633 & 26.693 \\
					\hline
					(15,21) & 0.1 & 0.556 & 1.737 & 0.702 & 1.982 & 1.046 & 2.484 \\
					& 0.3 & 7.820 & 12.486 & 8.430 & 12.970 & 9.605 & 13.915 \\
					& 0.5 & 16.712 & 21.803 & 17.310 & 22.134 & 18.444 & 22.788 \\
					& 0.7 & 23.525 & 27.100 & 23.951 & 27.275 & 24.796 & 27.639 \\
					& 0.9 & 27.427 & 28.581 & 27.712 & 28.779 & 28.329 & 29.124 \\
					\hline\hline
				\end{tabular}
			\end{table}

			\section{Data Analysis}\label{sec5}
			In this section, we present two real-life data analysis to illustrate the applicability of the derived results for the Lomax distribution. 
			
			\subsection{Oil gap-paper insulation failure data}
			We consider a datasets from \cite{jiang2023study}, which reports the breakdown voltage and failure time of oil-gap–paper insulation samples, measured over the same voltage duration using accelerated electrical–thermal ageing tests. Table~\ref{d1} presents two failure time datasets corresponding to oil gap--paper insulation specimens tested at different temperatures while maintaining the same voltage duration. Throughout this analysis, the voltage duration is assumed to be 20 seconds.
			
			\begin{table}[h!]
				\centering
				\caption{Failure time data of oil gap--paper insulation specimens at $80^\circ\mathrm{C}$ and $60^\circ\mathrm{C}$} \label{d1}
				\begin{tabular}{@{}ccc @{}}
					\Xhline{1pt}
					&T/$^{\circ}$C & Failure time/s \\
					\hhline{===}
					Data 1 & 80 & 283,237,211,227,278,254,242,264,273,294 \\
					Data 2 & 60 & 332,278,257,261,278,319,256,326,308,317 \\
					\Xhline{1pt}
				\end{tabular}
			\end{table}
			
			We used the parametric bootstrap Kolmogorov–Smirnov (KS) goodness-of-fit test to examine whether the datasets follow a multivariate Lomax distribution for various values of $b$. For each value of $b$, the corresponding $p$-values were obtained for Datasets 1 and 2. Based on these $p$-values, the data follow a two parameter exponential distribution. Table~\ref{d0} presents the $p$-values for Datasets 1 and 2 corresponding to the different values of $b$.
			
			\begin{table}[h!]
				\centering
				\caption{Bootstrap KS test $p$-values for different values of $b$.}
				\label{d0}
				\begin{tabular}{c|cc}
					\Xhline{1pt}
					$b$ & Dataset 1 & Dataset 2 \\
					\hhline{===}
					3 & 0.0561 & 0.0871 \\
					5 & 0.1341 & 0.1481 \\
					7 & 0.1632 & 0.1812 \\
					\Xhline{1pt}
				\end{tabular}
			\end{table}
			
			Let $X_1$ and $X_2$ denote the independent lifetimes of oil gap-paper insulation subjected to temperatures of $80^\circ\mathrm{C}$ and $60^\circ\mathrm{C}$, respectively. For given $\tau>0$, we have $X_1$ and $X_2$ be $\mbox{Exp}\left(\mu_1,\frac{\sigma_1}{\tau}\right)$ and $\mbox{Exp}\left(\mu_2,\frac{\sigma_2}{\tau}\right)$ respectively where $\mu_1$ and $\mu_2$ are location parameters. \cite{jiang2023study} reported that the failure time of oil gap-paper insulation decreases as the operating temperature increases. Consequently, the insulation operating at $80^\circ\mathrm{C}$ is expected to have a higher failure rate than that operating at $60^\circ\mathrm{C}$. Therefore, we consider 
			$\frac{\tau}{\sigma_1}\geq\frac{\tau}{\sigma_2}$, that is $\sigma_1\le\sigma_2$.

			Based on these datasets, we have obtained the values of the estimators for $\sigma_1$ and $\sigma_2$ of the multivariate Lomax distribution, as presented in Tables \ref{table11} and \ref{table12} with respect to the loss functions $L_1(\cdot)$ and $L_2(\cdot)$ respectively.

			\begin{table}[h!]
				\centering
				\caption{Estimated values of the improved estimators of $\sigma_1$ and $\sigma_2$ under $L_1(\cdot)$}
				\label{table11}
				\begin{tabular}{c|ccccc|ccc}
					\Xhline{1pt} 
					& \multicolumn{5}{c|}{Estimators of $\sigma_1$}
					& \multicolumn{3}{c}{Estimators of $\sigma_2$} \\ \cline{2-9}
					$b$ &
					$\delta_{11}$ &
					$\delta_{11}^1$ &
					$\delta_{12}^1$ &
					$\delta_{13}^1$ &
					$\delta_{\varphi_*^1}$ &
					$\delta_{12}$ &
					$\delta_{21}^1$ &
					$\delta_{\psi_*^1}$ \\
					\hhline{=========}
					3 & 130.770 & 115.523 & 130.770 & 130.770 & 101.918 & 107.387 & 115.523 & 140.998\\
					5 & 238.752 & 210.915 & 238.752 & 238.752 & 186.076 & 196.061 & 210.915 & 257.426\\
					7 & 345.984 & 305.645 & 345.984 & 345.984 & 269.649 & 284.120 & 305.645 & 373.046\\
					\Xhline{1pt} 
				\end{tabular}
			\end{table}
			
			\begin{table}[h!]
				\centering
				\caption{Estimated values of the improved estimators of $\sigma_1$ and $\sigma_2$ under $L_2(\cdot)$}
				\label{table12}
				\begin{tabular}{c|ccccc|ccc}
					\Xhline{1pt} 
					& \multicolumn{5}{c|}{Estimators of $\sigma_1$}
					& \multicolumn{3}{c}{Estimators of $\sigma_2$} \\ \cline{2-9}
					$b$ &
					$\delta_{21}$ &
					$\delta_{11}^2$ &
					$\delta_{12}^2$ &
					$\delta_{13}^2$ &
					$\delta_{\varphi_*^2}$ &
					$\delta_{22}$ &
					$\delta_{21}^2$ &
					$\delta_{\psi_*^2}$ \\
					\hhline{=========}
					3 & 100.667 &   91.667 &  100.667 &  100.667 &   80.546 &82.667 &   86.842 &  140.998 \\
					5 & 201.333 &  183.333 &  201.333 &  201.333 &  161.092  &165.333 &  173.684 &  257.426 \\
					7 & 302 &  275 &  302&  302 & 241.638 &248 &  260.526 & 373.046 \\
					\Xhline{1pt} 
				\end{tabular}
			\end{table}
			
			\subsection{Analysis of business failure data}
			In this example, we consider a business failure data set. The data were originally compiled by \cite{hutchinson1938study}, who studied the survival of different categories of business enterprises established in Poughkeepsie, New York, during the period 1844–1926. Later, \cite{lomax1954business} analyzed these data to investigate business mortality and observed that, unlike many engineering failure processes, the probability of business failure decreases as the age of the business increases.
			
			For our analysis, we consider two categories of businesses, namely retailers and craft. The cumulative probabilities of business failure over the first ten years of operation are presented in Table \ref{d2}. We used parametric bootstrap KS goodness-of-fit test to examine whether the datasets follow a $\mbox{Exp}\left(\mu,\frac{\sigma}{\tau}\right)$ for a given $\tau>0$.
			The test was performed for fixed values $b$, such as $b=3,5$ and $7$. The resulting bootstrap $p$-values are reported in Table \ref{d00}. Since all the obtained $p$-values exceed the $5\%$ significance level, there is no statistical evidence against the proposed model.
			
			\begin{table}[h!]
				\centering
				\caption{Probability density of business of failure in Poughkeepsie between 1844 and 1926} \label{d2}
				\begin{tabular}{@{}ccc @{}}
					\Xhline{1pt}
					Age in years& Retail (Dataset 1) & Craft (Dataset 2)  \\
					\hhline{===}
					0 & 0.57 & 0.5\\
					1 & 0.175 & 0.213 \\
					2 & 0.107 & 0.102 \\
					3 & 0.071 & 0.072 \\
					4 & 0.056 & 0.056 \\
					5 & 0.045 & 0.044 \\
					6 & 0.036 & 0.034 \\
					7 & 0.028 & 0.028 \\
					8 & 0.023 & 0.024 \\
					9 & 0.020 & 0.020 \\
					10 & 0.018 & 0.017 \\
					\Xhline{1pt}
				\end{tabular}
			\end{table}
			
			\begin{table}[h!]
				\centering
				\caption{Bootstrap KS test $p$-values for different values of $b$}
				\label{d00}
				\begin{tabular}{c|cc}
					\Xhline{1pt}
					$b$ & Dataset 1 & Dataset 2 \\
					\hhline{===}
					3 & 0.4454 & 0.5075 \\
					5 & 0.2763 & 0.2713 \\
					7 & 0.2002 & 0.2272 \\
					\Xhline{1pt}
				\end{tabular}
			\end{table}
			
			Let $X_1$ and $X_2$ denote the lifetimes of retail and craft business respectively. For given $\tau>0$, $X_1$ and $X_2$ be $\mbox{Exp}\left(\mu_1,\frac{\sigma_1}{\tau}\right)$ and $\mbox{Exp}\left(\mu_2,\frac{\sigma_2}{\tau}\right)$ respectively where $\mu_1$ and $\mu_2$ are location parameters. The probability densities of business failure in Poughkeepsie during the period $1844-1926$ show that craft businesses consistently exhibit a higher probability of failure than retail businesses at every business age from 1 to 10 years. This indicates that craft businesses are expected to have a higher failure rate than retail businesses. Therefore, $\frac{\tau}{\sigma_1}\geq\frac{\tau}{\sigma_2}$ so we have $\sigma_1\leq\sigma_2$.
			
			Based on these datasets, we have obtained the values of the estimators for $\sigma_1$ and $\sigma_2$ of the multivariate Lomax distribution, as presented in Tables \ref{ta11} and \ref{ta12} under the loss function $L_1(\cdot)$ and $L_2(\cdot)$ respectively.

			\begin{table}[h!]
				\centering
				\caption{Estimated values of the improved estimators of $\sigma_1$ and $\sigma_2$ under $L_1(\cdot)$}
				\label{ta11}
				\begin{tabular}{c|ccccc|ccc}
					\Xhline{1pt} 
					& \multicolumn{5}{c|}{Estimators of $\sigma_1$} & \multicolumn{3}{c}{Estimators of $\sigma_2$} \\ \cline{2-9}
					$b$ &
					$\delta_{11}$ &
					$\delta_{11}^1$ &
					$\delta_{12}^1$ &
					$\delta_{13}^1$ &
					$\delta_{\varphi_*^1}$ &
					$\delta_{12}$ &
					$\delta_{21}^1$ &
					$\delta_{\psi_*^1}$ \\
					\hhline{=========}
					3 &0.246&0.235&0.246&0.246&0.203& 0.238&0.238&0.289\\
					5 &0.448&0.43&0.448&0.448&0.37& 0.435&0.435&0.528\\
					7 &0.65&0.623&0.65&0.65&0.536& 0.631&0.631&0.765\\
					\Xhline{1pt} 
				\end{tabular}
			\end{table}
			
			\begin{table}[h!]
				\centering
				\caption{Estimated values of the improved estimators of $\sigma_1$ and $\sigma_2$ under $L_2(\cdot)$}
				\label{ta12}
				\begin{tabular}{c|ccccc|ccc}
					\Xhline{1pt} 
					& \multicolumn{5}{c|}{Estimators of $\sigma_1$}
					& \multicolumn{3}{c}{Estimators of $\sigma_2$} \\ \cline{2-9}
					$b$ &
					$\delta_{21}$ &
					$\delta_{11}^2$ &
					$\delta_{12}^2$ &
					$\delta_{13}^2$ &
					$\delta_{\varphi_*^2}$ &
					$\delta_{22}$ &
					$\delta_{21}^2$ &
					$\delta_{\psi_*^2}$ \\
					\hhline{=========}
					3 &0.19&0.187& 0.19&0.19&0.16 &0.185&0.185&0.289 \\
					5 &0.38&0.375&0.38&0.38&0.321&0.369&0.369&0.528\\
					7 &0.571&0.562&0.571&0.571&0.481&0.554&0.554&0.765\\
					\Xhline{1pt} 
				\end{tabular}
			\end{table}
			
			\section{Concluding remarks}	 
			In this work, we have investigated the problem of estimating the ordered scale parameters $\sigma_1$ and $\sigma_2$ of two scale mixture of exponential distributions under two scale invariant loss functions, namely the Stein loss and the symmetric loss. We consider a class of scale invariant estimators and propose sufficient conditions to obtain estimators that dominate the BAEE. Furthermore, we derive several Stein-type estimators that utilize all the information contained in samples from both populations and demonstrate that they outperform the BAEE. 
			
			Next, we have obtained a class of improved estimators. It has been observed that the boundary estimator of this class is a generalized Bayes estimator and it coincides with the Brewster and Zidek-type estimator. As an application, we have obtained improved estimators of the ordered scale parameter of two multivariate Lomax distributions and two exponential inverse Gaussian distributions.
			
			A numerical study was conducted to compare the risk performance of the proposed estimators. Furthermore, it was observed that, under the loss function $L_2(\cdot)$, the estimator $\delta_{\varphi_*^2}$ outperforms $\delta_{11}^2$ for $0<\eta<0.5$. However, when $\eta>0.5$, the estimator $\delta_{11}^2$ performs better than $\delta_{\varphi_*^2}$. Also, for estimating $\sigma_2$, the numerical results indicate that the estimator $\delta_{\psi^2}$ consistently outperforms $\delta_{21}^2$. Hence, $\delta_{\psi^2}$ is preferable to $\delta_{21}^2$ for estimating $\sigma_2$. A similar pattern of risk performance is observed for the exponential inverse Gaussian distribution. Finally to demonstrate the practical applicability of the proposed improved estimators, we consider two real data examples.
			\section{Appendix}
			We will use the following lemma to prove our results. For the sake of completeness, we state it below.
			\begin{lemma} [\cite{bobotas2009strawderman}]\label{lmmap}
				Let $f_1(x)$ and $f_2(x)$ be the densities supported on the domain set $\Omega_1$ and $\Omega_2$ respectively, where $\Omega_1 \subset \Omega_2$ and the ration $f_1(x)/f_2(x)$ is nondecreasing in $x \in \Omega_1$. If $X$ is a random variable having density $f_1(x)$ or $f_2(x)$ and also $a(x)$, $x \in \Omega_1$, is nondecreasing (nonincreasing) then $E_{f_1}a(X)\ge (\le)E_{f_2}a(X)$. Moreover, if $a(x)$ and the ration $f_1(x)/f_2(x)$ are strictly monotone, then $E_{f_1}a(X)>(<)E_{f_2}a(X)$.
			\end{lemma}
			
			\subsection{Proof of Theorem \ref{thIERD}}\label{proof1}
			\begin{enumerate}
				\item [(i)] The risk difference of the estimators $\delta_{11}$ and $\delta_{\varphi_1}$ is 
				\begin{align*}
					\Delta_1 =& E\left[\int_{1}^{\infty}\left(\varphi_1'\left(\frac{1}{\eta}\frac{v_2}{v_1}t\right)\frac{1}{\eta}v_2-\frac{\varphi_1'\left(\frac{1}{\eta}\frac{v_2}{v_1}t\right)\frac{1}{\eta}\frac{v_2}{v_1}}{\varphi_1^2\left(\frac{1}{\eta}\frac{v_2}{v_1}t\right)v_1}\right)dt\right]\\
					\geq&E\left[\int_{1}^{\infty}\left(v_2-\frac{v_2}{v_1^2\varphi_1^2\left(t\frac{v_2}{v_1}\right)}\right)\varphi_1'\left(\frac{1}{\eta}\frac{v_2}{v_1}t\right)dt\right] 
					\\
					=&\int_{0}^{\infty}\int_{0}^{\infty}\int_{0}^{\infty}\int_{1}^{\infty}\left(v_2-\frac{v_2}{v_1^2\varphi_1^2\left(t\frac{v_2}{v_1}\right)}\right)\varphi_1'\left(\frac{1}{\eta}\frac{v_2}{v_1}t\right)\tau g_1(\tau v_1) \tau g_2(\tau v_2)dt dv_1 dv_2dH(\tau)
				\end{align*}		
				Making the transformation $u=t\frac{v_2}{v_1}$ and $x=\frac{uv_1}{t}$, the we get
				\begin{align*}
					\Delta_1 \geq&\int_{0}^{\infty}\int_{0}^{\infty}\int_{0}^{\infty}\int_{0}^{uv_1} \left(1-\frac{1}{v_1^2\varphi_1^2\left(u\right)}\right)v_1\varphi_1'\left(\frac{u}{\eta}\right)\tau g_1(\tau v_1) \tau g_2(\tau x)dx dv_1 du dH(\tau)\\
					=&\int_{0}^{\infty}\varphi_1'\left(\frac{u}{\eta}\right)\int_{0}^{\infty}\int_{0}^{\infty} \left(1-\frac{1}{v_1^2\varphi_1^2\left(u\right)}\right)v_1\tau g_1(\tau v_1)\int_{1}^{\tau uv_1} g_2(x) dx dv_1 dH(\tau) du
				\end{align*}		
				Now the risk difference $\Delta_1\geq0$ if 
				\begin{align*}
					\varphi_1(u)\geq \left( \frac{\int_{0}^{\infty}\int_{0}^{\infty}\frac{1}{v_1}\tau g_1(\tau v_1)G_2(\tau uv_1)dv_1 dH(\tau)} {\int_{0}^{\infty}\int_{0}^{\infty}v_1\tau g_1(\tau v_1)G_2(\tau uv_1)dv_1 dH(\tau)}\right)^{1/2}
				\end{align*}
				Put $\tau v_1= y$, we have
				\begin{align*}
					\varphi_1(u)\geq \left( \frac{\int_{0}^{\infty}\tau \int_{0}^{\infty}\frac{1}{y} g_1(y)G_2(yu)dy dH(\tau)} {\int_{0}^{\infty}\frac{1}{\tau}\int_{0}^{\infty}y g_1(y)G_2(yu)dy dH(\tau)}\right)^{1/2}= \left(\frac{E(\tau)\int_{0}^{\infty} y^{p_1-3}e^{-y} \int_{0}^{yu}x^{p_2-2}e^{-x}dx dy}{E(1/\tau)\int_{0}^{\infty} y^{p_1-1}e^{-y} \int_{0}^{yu}x^{p_2-2}e^{-x}dx dy}\right)^{1/2}
				\end{align*}
				
				Take $z=\frac{x}{yu}$, we have 
				\begin{align*}
					\varphi_1(u)\geq & \left(\frac{E(\tau)\int_{0}^{\infty} y^{p_1+p_2-4} \int_{0}^{1}z^{p_2-2}e^{-y(1+zu)}dz dy}{E(1/\tau)\int_{0}^{\infty} y^{p_1+p_2-2} \int_{0}^{1}z^{p_2-2}e^{-y(1+zu)}dz dy}\right)^{1/2}\\
				%
					=& \left(\frac{E(\tau) B\left(\frac{u}{1+u}; p_1-2,p_2-1\right)}{E(1/\tau)(p_1+p_2-2)(p_1+p_2-3)B\left(\frac{u}{1+u}; p_1,p_2-1\right)}\right)^{1/2}=\varphi_{*}^1(u).
				\end{align*}
				
				\item [(ii)] The proof for the Stein loss function $L_2(\cdot)$ is similar to (i). So we omit it. 
			\end{enumerate}
			
			\subsection{Proof of Theorem \ref{ths21}}\label{proof2}
			
			\begin{enumerate}
				\item [(i)]  Under the loss function $L_1(\cdot)$, risk of the estimator $\delta_{\varphi_2}(W,W_1)$ can be expressed as 
				\begin{align*}
					R(\delta_{\varphi_2},\mu_1,\sigma_1,\sigma_2)=E^{W,W_1}E\left[\left(\varphi_2(W,W_1)V_1+\frac{1}{ \varphi_2(W,W_1)V_1} - 2\right)\bigg| W, W_1\right]
				\end{align*}
				For given $\tau>0$, the conditional density of $V_1$ given $W=w$, $W_1=w_1$ obtain as
				\begin{equation*}
					f_{\eta,\rho}(v_1|w,w_1) = \frac{\int_{0}^{\infty}v_1^{p_1+p_2-2}e^{-\tau v_1(1+w\eta+p_1w_1)} e^{p_1\tau\rho}\tau^{p_1+p_2-1}dH(\tau)} {\int_{0}^{\infty}\int_{\max\{0,\frac{\rho}{w_1}\}}^{\infty}v_1^{p_1+p_2-2}e^{-\tau v_1(1+w\eta+p_1w_1)} e^{p_1\tau\rho}\tau^{p_1+p_2-1}dv_1 dH(\tau)} 
				\end{equation*}
				where $\eta$=$\frac{\sigma_1}{\sigma_2}\leq1$, $\rho=\frac{\mu_1}{\sigma_1}\in \mathbb{R}$. It can be easily seen that the conditional risk function 
				\begin{equation*}
					R_1(\delta_{\varphi_2},\eta,\rho)=	E\left[\left(\varphi_2(w,w_1)V_1+\frac{1}{\varphi_2(w,w_1)V_1} - 2\right)\bigg|W=w,W_1=w_1\right]
				\end{equation*}
				is minimized at 
				\begin{equation*}
					\varphi_2(w,w_1;\eta,\rho)= \left(\frac{E\left[ 1/V_1|W=w, W_1=w_1\right]}{E\left[V_1|W=w, W_1=w_1\right]}\right)^{1/2}
				\end{equation*}
				Now we will consider two cases. First case, we consider $\mu_1>0$, $w_1>0$. In this case we have
				\begin{equation*}
					\varphi_2(w,w_1;\eta,\rho) =\left( \frac{\int_{0}^{\infty}\int_{\frac{\rho}{w_1}}^{\infty}v_1^{p_1+p_2-3}e^{-\tau v_1(1+w\eta+p_1w_1)} e^{p_1\tau \rho} \tau^{p_1+p_2-1}dv_1dH(\tau)}{\int_{0}^{\infty}\int_{\frac{\rho}{w_1}}^{\infty}v_1^{p_1+p_2-1}e^{-\tau v_1(1+w\eta+p_1w_1)} e^{p_1\tau \rho} \tau^{p_1+p_2-1}dv_1dH(\tau)}\right)^{1/2}.
				\end{equation*}
				Now by taking the transformation $z_1=\tau v_1(1+w\eta+p_1w_1)$, we obtain
				\begin{equation*}
					\varphi_2(w,w_1;\eta,\rho) =(1+w\eta+p_1w_1)\left( \frac{\int_{0}^{\infty}\tau e^{p_1\tau \rho }\int_{\xi}^{\infty}z_1^{p_1+p_2-3}e^{-z_1}dz_1dH(\tau)}{\int_{0}^{\infty}1/\tau e^{p_1\tau \rho}\int_{\xi}^{\infty}z_1^{p_1+p_2-1}e^{-z_1} dz_1dH(\tau)}\right)^{1/2}.
				\end{equation*}
				For given $\tau>0$, we can easily seen that
				\begin{equation*}
					\frac{\int_{\xi}^{\infty}z_1^{p_1+p_2-3}e^{-z_1}dz_1}{\int_{\xi}^{\infty}z_1^{p_1+p_2-1}e^{-z_1} dz_1}=E_{\xi}(Z_1^{-2}),
				\end{equation*}
				where $Z_1$ has density $g(z_1,\xi)\propto z_1^{p_1+p_2-1}e^{-z_1}I_{(\xi,\infty)}(z_1)$ with $\xi= \frac{\tau \rho}{w_1}(1+w\eta+p_1w_1)>0$.
				For $\xi>0$, $\frac{g(z_1,\xi)}{g(z_1,0)}$ is non-decreasing then we have $E_{\xi}Z_1^{-2}\leq E_0Z_1^{-2}=\frac{1}{(p_1+p_2-1)(p_1+p_2-2)}$ and hence we get  $\int_{\xi}^{\infty}z_1^{p_1+p_2-3}e^{-z_1}dz_1\leq \frac{1}{(p_1+p_2-1)(p_1+p_2-2)}\int_{\xi}^{\infty}z_1^{p_1+p_2-1}e^{-z_1}dz_1$.
				Now we obtain that
				\begin{align*}
					\varphi_2(w,w_1;\eta,\rho) &\leq\frac{(1+w\eta+p_1w_1)}{\sqrt{(p_1+p_2-1)(p_1+p_2-2)}}\left( \frac{\int_{0}^{\infty}\tau e^{p_1\tau\rho}\int_{\xi}^{\infty}z_1^{p_1+p_2-1}e^{-z_1}dz_1dH(\tau)}{\int_{0}^{\infty}\frac{1}{\tau}e^{p_1\tau\rho}\int_{\xi}^{\infty}z_1^{p_1+p_2-1}e^{-z_1} dz_1dH(\tau)}\right)^{1/2}
				\end{align*}
				Again we take the transformation $z_1=\tau x$ then we have 
				\begin{align} \nonumber \label{eq3.8}
					\varphi_2(w,w_1;\eta,\rho) &\leq \frac{(1+w\eta+p_1w_1)}{\sqrt{(p_1+p_2-1)(p_1+p_2-2)}}\left( \frac{\int_{0}^{\infty}\int_{\xi_1}^{\infty}x^{p_1+p_2-1}\tau^{p_1+p_2+1} e^{p_1\tau\rho}e^{-\tau x}dxdH(\tau)}{\int_{0}^{\infty}\int_{\xi_1}^{\infty}x^{p_1+p_2-1}\tau^{p_1+p_2-1} e^{p_1\tau\rho}e^{-\tau x} dxdH(\tau)}\right)^{1/2}\\
					&= \frac{(1+w\eta+p_1w_1)}{\sqrt{(p_1+p_2-1)(p_1+p_2-2)}}\left( \frac{\int_{0}^{\infty}\int_{\xi_1}^{\infty}x^{p_1+p_2-1}\tau^{p_1+p_2+1} e^{-(x-p_1\rho)\tau}dxdH(\tau)}{\int_{0}^{\infty}\int_{\xi_1}^{\infty}x^{p_1+p_2-1}\tau^{p_1+p_2-1} e^{-(x-p_1\rho)\tau} dxdH(\tau)}\right)^{1/2}
				\end{align}
				where, $\xi_1=\frac{\rho}{w_1}(1+w\eta+p_1w_1)>0$. Now for $k=x-p_1\rho>0$, we obtain
				\begin{align*}
					\frac{\int_{0}^{\infty}\tau^{p_1+p_2+1} e^{-k\tau}dH(\tau)}{\int_{0}^{\infty}\tau^{p_1+p_2-1} e^{-k\tau}dH(\tau)}=\int_{0}^{\infty}\tau^2f_k(\tau)dH(\tau)
				\end{align*}
				where $f_k(\tau) \propto \tau^{p_1+p_2-1} e^{-k\tau} $. But, $\frac{f_k(\tau)}{f_0(\tau)}$ is decreasing in $\tau$, so by using Lemma (\ref{lmmap}) we have, 
				\begin{align}\label{eq3.9}
					\int_{0}^{\infty}\tau^2f_k(\tau)dH(\tau) \leq \int_{0}^{\infty}\tau^2f_0(\tau)dH(\tau)=\frac{E(\tau^{p_1+p_2+1})}{E(\tau^{p_1+p_2-1})}
				\end{align} 
				From (\ref{eq3.8}) and (\ref{eq3.9}) with $\eta\leq 1$ we get
				\begin{align}\label{eq3.10}
					\varphi_2(w,w_1;\eta,\rho) &\leq \frac{(1+w+p_1w_1)}{\sqrt{(p_1+p_2-1)(p_1+p_2-2)}}\left(\frac{E(\tau^{p_1+p_2+1})}{E(\tau^{p_1+p_2-1})}\right)^{1/2}
				\end{align}
				
				\noindent Now we consider the case $\mu_1 \le 0$, $w_1>0$. In this case we have, 
				$$\varphi_2(w,w_1;\eta,\rho) =\left( \frac{\int_{0}^{\infty}\int_{0}^{\infty}v_1^{p_1+p_2-3}e^{-\tau v_1(1+w\eta+p_1w_1)} e^{\rho p_1\tau} \tau^{p_1+p_2-1}dv_1dH(\tau)}{\int_{0}^{\infty}\int_{0}^{\infty}v_1^{p_1+p_2-1}e^{-\tau v_1(1+w\eta+p_1w_1)} e^{\rho p_1\tau} \tau^{p_1+p_2-1}dv_1dH(\tau)}\right)^{1/2}$$
				After using the transformation $z_1=\tau v_1(1+w\eta+p_1w_1)$ we have,
				\begin{align}\label{eq3.11}\nonumber
					\varphi_2(w,w_1;\eta,\rho) =&(1+w\eta+p_1w_1)\left( \frac{\int_{0}^{\infty}\tau e^{p_1\rho \tau}\int_{0}^{\infty}z^{p_1+p_2-3}e^{-z}dzdH(\tau)}{\int_{0}^{\infty}1/\tau e^{p_1\tau \rho}\int_{0}^{\infty}z^{p_1+p_2-1}e^{-z} dzdH(\tau)}\right)^{1/2}\\\nonumber
					=& \frac{(1+w\eta+p_1w_1)}{\sqrt{(p_1+p_2-1)(p_1+p_2-2)}}\left( \frac{\int_{0}^{\infty}\tau e^{p_1\rho \tau}dH(\tau)}{\int_{0}^{\infty}1/\tau e^{p_1\tau \rho} dH(\tau)}\right)^{1/2}\\
					=&\frac{(1+w\eta+p_1w_1)}{ \sqrt{(p_1+p_2-1)(p_1+p_2-2)}}\left(\int_{0}^{\infty}\tau^2 f_{\rho}(\tau)\right)^{1/2}
				\end{align}
				where $f_{\rho}(\tau)=\frac{\frac{1}{\tau}e^{n\tau \rho}}{\int_{0}^{\infty}\frac{1}{\tau}e^{\tau n\rho}dH(\tau)}$.
				Set $f_0(\tau)=\frac{\tau^{-1}}{\int_{0}^{\infty}\tau^{-1}dH(\tau)}$. Now $\frac{f_{\rho}(\tau)}{f_0(\tau)}$ is decreasing in $\tau$, then by Lemma (\ref{lmmap}) we have 
				\begin{align}\label{eq3.12}
					\int_{0}^{\infty}\tau^2f_{\rho}(\tau)dH(\tau) \leq \int_{0}^{\infty}\tau^2f_0(\tau)dH(\tau)=\frac{E(\tau)}{E(\tau^{-1})}.
				\end{align} 
				Hence from equation (\ref{eq3.11}) and (\ref{eq3.12}) we obtain
				\begin{align}\label{eq3.13}
					\varphi_2(w,w_1;\eta,\rho) &\leq \frac{(1+w+p_1w_1)}{\sqrt{(p_1+p_2-1)(p_1+p_2-2)}}\left(\frac{E(\tau)}{E(\tau^{-1})}\right)^{1/2}.
				\end{align}
				
				Hence from the equations (\ref{eq3.13}) and (\ref{eq3.10}) we get,   
				\begin{align*}
					\varphi_2(w, w_1;\eta,\rho) \leq \frac{\big(1 + w + p_1 w_1\big)\min\Bigg\{ 
						\left( \frac{E(\tau^{p_1 + p_2 + 1})}
						{E(\tau^{p_1 + p_2 - 1})} \right)^{1/2}, \left( \frac{E(\tau)}
						{E(\tau^{-1})} \right)^{1/2}
						\Bigg\}}{\sqrt{(p_1+p_2-1)(p_1+p_2-2)}} =\varphi_{12}(w,w_1).
				\end{align*}
				Now using the convexity of R$_1(\delta_{\varphi_2},\eta,\rho)$ for $P(\varphi_{12}(W,W_1)<\varphi_2(W,W_1))>0$ we get the result. 
				%
				\item [(ii)] The proof for the $L_2(\cdot)$ loss function is similar to $(i)$. So we omit it.  
			\end{enumerate}
			
			\section*{Data availability statement}
			Data used in this study are publicly available in the literature and cited in this manuscript.
			\section*{Disclosure statement}
			No potential conflict of interest was reported by the authors.
			\section*{Funding}
			Lakshmi Kanta Patra thanks the Anusandhan National Research Foundation, India for providing financial support to carry out this research with project number MTR/2023/000229.

		\end{document}